\documentclass[preprint,12pt]{elsarticle}

\usepackage{multirow}    

\usepackage{amsmath}
\usepackage{amssymb}

\usepackage{color}
\usepackage{array}
\usepackage{booktabs}
\usepackage{enumitem}
\usepackage[misc,geometry]{ifsym}

\DeclareMathAlphabet\mathbfcal{OMS}{cmsy}{b}{n} 

\usepackage{graphicx}
\usepackage{mathrsfs}
\usepackage{color}
\usepackage[bookmarks,bookmarksnumbered]{hyperref}
\hypersetup{colorlinks = true,
	linkcolor = blue,
	anchorcolor =red,
	citecolor = blue,
	filecolor = red,
	urlcolor = red,
	pdfauthor=author}

\usepackage{tikz}
\usetikzlibrary{calc,positioning}

\journal{Mechanical Systems and Signal Processing}

\begin{document}

\begin{frontmatter}



\title{Parallelized computation of quasi-periodic solutions for finite element problems: A Fourier series expansion-based shooting method}


\author[inst1]{Junqing Wu} 
\author[inst1]{Ling Hong} 
\author[inst2]{Mingwu Li\corref{cor1}}
\ead{limw@sustech.edu.cn}
\author[inst1]{Jun Jiang\corref{cor1}} 
\ead{jun.jiang@xjtu.edu.cn}

\affiliation[inst1]{organization={State Key Laboratory for Strength and Vibration, Xi’an Jiaotong University},
            addressline={No.28, West Xianning Road}, 
            city={Xi’an},
            postcode={710049}, 
            country={China}}
            
\affiliation[inst2]{organization={Department of Mechanics and Aerospace Engineering, Southern University of Science and Technology},
	addressline={No.1088 Xueyuan Avenue}, 
	city={Shenzhen},
	postcode={518055}, 
	country={China}}

\cortext[cor1]{Corresponding author}
	
\begin{abstract}
High-dimensional nonlinear mechanical systems admit quasi-periodic solutions that are essential for the understanding of dynamical systems. These quasi-periodic solutions stay on some invariant tori governed by complex PDEs in hyper-time. Here, we propose a Fourier series expansion-based shooting method (FSE-Shooting) for the parallelized computation of quasi-periodic solution with $d$ base frequencies ($d \ge 2$). We represent the associated $d$-torus as a collection of trajectories initialized at a ($d-1$)-torus. We derive a set of ODEs that hold for any of these trajectories. We also derive a set of boundary conditions that couple the initial and terminal states of these trajectories and then formulate a set of nonlinear algebraic equations via the coupling conditions. We use Fourier series expansion to parameterize the ($d-1$)-torus and shooting method to iterate the Fourier coefficients associated with initial torus such that the coupling conditions are satisfied. In particular, the terminal points of these trajectories are computed in parallel via Newmark integration, where the time points and Fourier coefficients are transformed to each other by alternating Frequency-Time method. A straightforward phase condition is devised to track the quasi-periodic solutions with a priori unknown base frequencies. Additionally, the by-product of the FSE-Shooting can also be directly used to compute the Lyapunov exponents to assess the stability of quasi-periodic solutions. The results of three finite element systems show the efficiency and versatility of FSE-Shooting in high-dimensional nonlinear dynamical systems, including a three-dimensional shell structure with $1872$ DOFs.

\end{abstract}


\end{frontmatter}



\section{Introduction}\label{sec1}

The determination of stationary solutions of nonlinear dynamical systems governed by the ordinary differential equations (ODEs), as well as the assessment of their stability, plays an important role in science and engineering \cite{R1,R2,R3}. The quasi-periodic solutions have been widely reported and investigated in the literature \cite{R4,R5,R6,R7,R8,R9,R10}, where engineering problems are often modeled as low-dimensional ODEs for computational convenience. This is because the dynamical behavior of quasi-periodic solutions is far more complex than that of equilibrium and periodic solutions. Unlike numerical methods for periodic solutions, which have been widely applied to high-dimensional systems \cite{R11,R12,R13,R14,R15,R16,R17}, the existing numerical methods for solving quasi-periodic solutions can only handle low-dimensional systems with two base frequencies.The few existing examples of solving quasi-periodic solutions for finite element problems are also based on their reduced-order models \cite{R18,R19,R20}. Therefore, this paper aims to develop an efficient numerical method for quasi-periodic solutions, imposing no restrictions on the number of base frequencies or whether they are known a priori. The proposed method is expected to offer theoretical support and methodological tools for exploring the dynamical behaviors of full models in dynamical systems.

Because of multiple incommensurate base frequencies, quasi-periodic solutions never repeat in the physical time domain. Such solutions are often investigated in the context of high-dimensional invariant tori in phase space \cite{R21,R22,R23}. For a solution with $d$ base frequencies, the parametrization of its $d$-torus \cite{R24,R25,R26,R27,R28} is key to numerical methods. Multi-period parametrization sets the derivative of each hyper-time variable with respect to physical time to one base frequency. In this case, the $d$-torus function satisfies partial differential equations (PDEs) and is subject to $d$ sets of periodic boundary conditions. There are two types of discretization techniques for this problem, namely full-discretization and semi-discretization methods.

 The former directly constructs a series of quasi-periodic shape functions with undetermined coefficients to approximate the $d$-torus, and formulates a set of nonlinear algebraic equations (NAEs) to iterate coefficients by imposing PDEs via the weighted-residual methods. For example, the multi-dimensional harmonic balance method (MHB) \cite{R29,R30,R31,R32,R33} define trigonometric functions with multi-time scales as shape functions and use Fourier–Galerkin projection to construct NAEs. By constructing Lagrange basis polynomials functions at the base points on the hyper-time domain, spectral collocation method (SCO) \cite{R34} requires that the PDEs are satisfied at the selected collocation-nodes, which also formulates the NAEs for unknowns of $d$-torus at based points. Although the finite difference method (FD) \cite{R23,R35,R36} involves no explicit shape functions, it constructs the NAEs similar to those in SCO by using finite differences at base points to approximate differentiations.

 Semi-discretization methods \cite{R25} provide an alternative discretization strategy, which uses periodic shape functions with variable coefficients to sequentially approximate the $d$-torus. By discretizing these periodic shape functions in the PDEs, one obtains new PDEs for variable coefficients, which are regarded as ($d−1$)-torus. After $d$ such steps, the variable coefficients are updated to constant coefficients, and the PDEs are eventually reformulated as NAEs. The conventional variable-coefficients harmonic balance method (VCHB) \cite{R37,R38,R39,R40} may be the typical examples for 2 base frequencies case. Recently, we proposed a multi-step variable-coefficients formulation (\textit{m}-VCF) \cite{R41} to establish a unified framework of semi-discretization method for solving quasi-periodic solutions with multiple base frequencies, where the discretization at each step can be carried out by means of any one of three typical discretization methods for periodic solutions, such as HB \cite{R42,R43,R44,R45}, CO \cite{R46,R47} and FD \cite{R48,R49}. Due to the decoupling of large matrix operations in the full-discretization methods, \textit{m}-VCF can efficiently evaluate the nonlinear terms in dynamical systems \cite{R41}. By the way, we also propose a single-step constant-coefficients formulation (\textit{s}-CCF) to extend the full-discretization method to the multi-frequency scenario in \ref{appA}.

However, the number of shape functions for $d$-torus increases exponentially with the number of base frequencies, which substantially increases the difficulty of computing the corresponding unknown coefficients, particularly for high-dimensional dynamical systems. Shooting-base method may be a good idea to reduce the number of unknowns in NAEs. For quasi-periodic solutions with two base frequencies, the extension of shooting idea can be traced back to the work of Kaas-Petersen \cite{R71,R72} and Kim \cite{R73}, whose NAEs aim to enforce that initial state of the quasi-periodic solution matches its higher-order Poincaré map. However, the initial state mapping can only be approximated via interpolation after long-time numerical integration.

More reliable extensions of the shooting method can be found in astrodynamics \cite{R74}, which focus on searching for the invariant curve of maps corresponding to 2-torus: one using surfaces of section \cite{R75}, and one using a stroboscopic map (also called the first-order Poincaré map) \cite{R76,R77}. A collection of trajectories initialized at invariant curve are used to parameterized 2-torus, which are governed by the original ODEs. Since the latter fixes the time of the trajectory return section, NAEs are constructed based on the fixed phase difference between their initial and final points, where the discrete Fourier transform (DFT) matrix and its inverse are employed. The same idea can also be found in \cite{R22,R50}, where each trajectory is discretized using the collocation method, which increases the number of unknowns in NAEs. An extended shooting method for $d$-torus is independently proposed in \cite{R35}, where points uniformly distributed on $(d-1)$-torus are defined as unknowns of NAEs. However, no explicit NAEs are derived, since cubic spline interpolation is used to fit $(d-1)$-torus.

In this work, we propose a Fourier series expansion-based shooting method (FSE-Shooting) for efficiently computing quasi-periodic solutions with arbitrary $d$ base frequencies. The multi-trajectory parametrization is defined to represent $d$-torus, whose hyper-time domain can be considered as a hyper-parallelepiped in that of previous multi-period parametrization. Since the cross-section of the $d$-torus preserves $d-1$ sets of periodic boundary conditions at any state along the trajectory, a Fourier series expansion with variable Fourier coefficients is used to approximate the $d$-torus. By means of the phase difference between initial and terminal states, we derive the coupling conditions between their Fourier coefficients. In order to define the Fourier coefficients of initial state as initial conditions for shooting process, the alternating Frequency-Time method (AFT) \cite{R39,R52,R53} is employed to convert between Fourier coefficients and their corresponding points on trajectories; the time integration method, such as Newmark integration \cite{R11,R51}, is used to compute terminal points from initial points along each trajectory.

It is also difficult to determine the quasi-periodic solution containing unknown base frequencies. Similar to the phase conditions of periodic solutions \cite{R56,R57,R58,R59}, several phase conditions are available to fix the phases of the unknown base frequencies for quasi-periodic solutions \cite{R39,R60,R61}. However, the existing phase conditions cannot guarantee robust solutions during one-parameter continuation \cite{R62,R63,R64,R65}, because they are heavily dependent on the accuracy of the predicted initial condition or the solution point of previous step. To overcome these drawbacks, by using the orthogonality in sense of the integral between current solution and its iterative variables at each step, we proposed a simple phase condition for robust continuation, which has been validated in both the MHB \cite{R66} and method constructed with \textit{m}-VCF \cite{R41}. In this work, the variant of this phase condition is also developed in FSE-Shooting to track the solutions with unknown base frequencies.

The main contribution of proposed FSE-Shooting can be concluded as follows: 1) It is applicable to any number of base frequencies, and the frequencies can be unknown; 2) It analytically derives the NAEs for Fourier coefficients, avoiding the frequency aliasing that may occur in previous shooting methods; 3) Minimizing the number of unknowns and designing parallel integration for trajectories effectively improve computational efficiency for high-dimensional systems.

The rest of this paper is organized as follows. Two types of parametrizations for $d$-torus and their intrinsic relationship are revealed in Section \ref{sec2}. Then, in Section \ref{sec3}, the Fourier series expansion-based shooting method is proposed based on a new set of boundary conditions of the Fourier variable-coefficients. A phase condition is also proposed to track the solution with unknown base frequencies. The feasibility of presented method is demonstrated by applying it to three nonlinear systems in Section \ref{sec4}. Finally, conclusions and future perspectives are drawn in the last section.

\section{Parametrizations of quasi-periodic solution and its dynamical systems}\label{sec2}

Consider the finite element model of a mechanical system, the governing equations of motion with $n$-DOFs are the following form

\begin{equation}
	\mathbf{M}\ddot{\mathbf{q}}(t)+\mathbf{D}\dot{\mathbf{q}}(t)+\mathbf{K}\mathbf{q}(t)+\mathbf{\Theta}[\mathbf{f}_{nl}(\mathbf{q},\dot{\mathbf{q}})-\mathbf{e}(\mathbf{\Omega}\times t)]=\mathbf{0},\label{2ODEs}
\end{equation}
where $t \in \mathbb{R}$ denotes physical time. $\mathbf{M}$, $\mathbf{D}$, $\mathbf{K} \in \mathbb{R}^{n \times n}$ are the generalized mass, the viscous damping and the stiffness matrices, each of size $n \times n$. The vector $\mathbf{f}_{nl}(\mathbf{q}, \dot{\mathbf{q}}) \in \mathbb{R}^{n \times 1}$ represents the nonlinear forces that depend on generalized displacements and generalized velocities. $\mathbf{e}$ is the vector of external excitation, with $\mathbf{\Omega} = [\omega_1, \dots, \omega_e]^\top \in \mathbb{R}^{e \times 1}$ being the $e$-dimensional excited frequencies. $\mathbf{\Theta}$ is the force distribution matrix, typically an $n \times n$ identity matrix. Over-dot denotes the first-order derivative with respect to physical time $t$.

Let $\mathbf{x}=[\mathbf{q};\dot{\mathbf{q}}]$, the second-order ODEs above can be rewritten as the first-order form below, with $2n$ state variables

\begin{equation}
	\dot{\mathbf{x}}(t)=\mathbf{f}\left(\mathbf{x}(t),\mathbf{\Omega}\times t\right),\quad\mathbf{f}:\mathbb{R}^{2n\times1}\times\mathbb{R}^{e\times1}\times\mathbb{R}\mapsto\mathbb{R}^{2n\times1},\label{1ODEs}
\end{equation}
Assume that system \eqref{1ODEs} has a family of quasi-periodic solutions $\mathbf{x}(t) \in \mathbb{R}^{2n \times 1}$, with $d$ base frequencies given by $\boldsymbol{\omega} = [\mathbf{\Omega}^\top, \omega_{e+1}, \dots, \omega_d]^\top \in \mathbb{R}^{d \times 1}$ and $d \geq e$, where $\boldsymbol{\omega}_u = [\omega_{e+1}, \dots, \omega_d]^\top$ are a priori unknown intrinsic frequencies.

\begin{figure}[ht] 
	\centering 
	\includegraphics[width=0.75\textwidth]{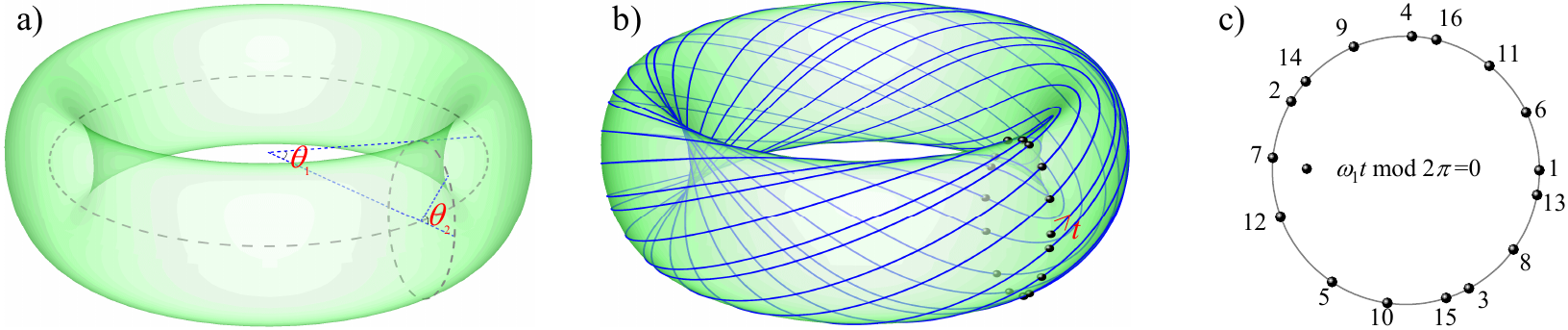} 
	\caption{\small A quasi-periodic solution with 2 base frequencies: a) A 2-torus; b) The trajectories of the quasi-periodic densely filling the 2-torus; c) The mapping points on the first-order Poincaré section (1-torus).} 
	\label{fig: 2-torus}
\end{figure}

The trajectory of a quasi-periodic solution will never repeat in the physical time domain but densely fill a $d$-dimensional invariant torus ($d$-torus). By setting torus coordinate $\boldsymbol{\theta} = [\theta_1, \dots, \theta_d]^\top$, each component is $\theta_i = \omega_i t \mod 2\pi$ for $i = 1, \dots, d$. Fig. \ref{fig: 2-torus} illustrates a quasi-periodic solution with 2 base frequencies that is a 2-torus with the trajectories filling densely on its torus surface. In Fig. \ref{fig: 2-torus}c, the first-order Poincaré section at $\theta_1 = 0$ of the torus, a 1-torus, is also drawn. For the general case, the first-order Poincaré section at $\theta_1 = 0$ of quasi-periodic solution with $d$ base frequencies is a $(d-1)$-torus associated with the torus coordinates $\tilde{\boldsymbol{\theta}} = [\theta_2, \dots, \theta_d]^\top$, on which two adjacent map points have the phase shift of $\tilde{\boldsymbol{\theta}} = 2\pi\boldsymbol{\rho}$, with $\boldsymbol{\rho} = [\rho_2, \dots, \rho_d]^\top$ and $\rho_j = \omega_j / \omega_1$, $j = 2, \dots, d$. In the following it will be shown that the way of parametrization on $d$-torus is crucial to effectively determine the quasi-periodic solution.

\subsection{Multi-period parametrization}\label{sec2.1}

\begin{figure}[ht] 
	\centering 
	\includegraphics[width=0.75\textwidth]{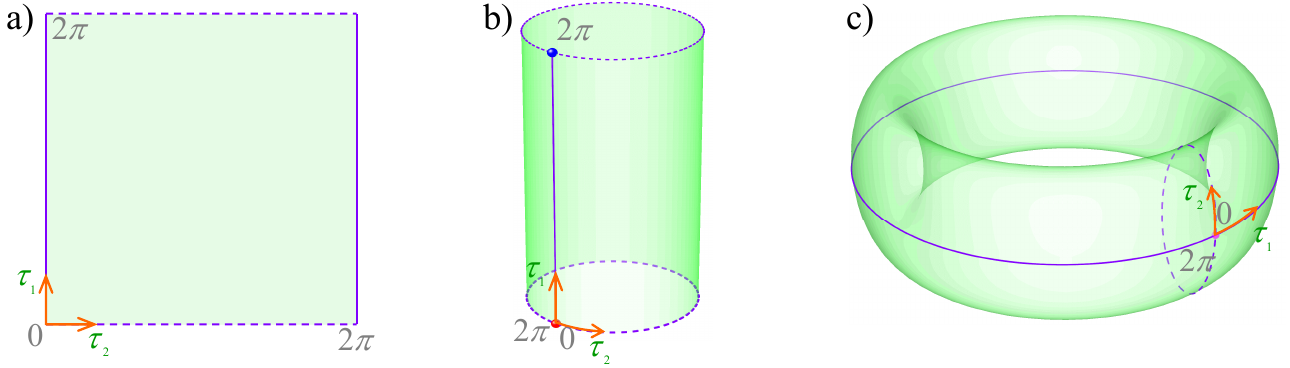} 
	\caption{\small The hyper-time domain $\mathbb{T}^d$ with $d = 2$: a) Depiction by an orthogonal basis; b) Connecting the pair of periodic boundaries of $\tau_2$; c) Further connecting the pair of periodic boundaries of $\tau_1$.} 
	\label{fig: tau_time}
\end{figure}

A popular parametrization of $d$-torus is based on the concept of hyper-time domain $\boldsymbol{\tau}\in\mathbb{T}^d=[0,2\pi)^d$, its hyper-time variables are required to be orthogonal to each other and satisfy the rule

\begin{equation}
	\frac{\mathrm{d}\tau_i}{\mathrm{d}t}=\omega_i,\,i=1,\cdots,d,\label{EqTau}
\end{equation}
For simplicity and without loss of generality, let $d$ hyper-time variables be $\tau_i = \omega_i t$, $i = 1, \dots, d$. So, the $d$-torus is typically parametrized by $\mathbf{y}(\boldsymbol{\tau}): \mathbb{T}^d \mapsto \mathbb{R}^{2n \times 1}$, satisfying $d$ periodic boundary conditions $\mathbf{y}(\dots, 0, \dots) = \mathbf{y}(\dots, 2\pi, \dots)$. For an intuitive understanding, Fig. \ref{fig: tau_time} exemplifies the equivalence between a 2-torus of the quasi-periodic solution with two base frequencies and the hyper-time domain with two orthogonal variables. Fig. \ref{fig: tau_time}a shows the hyper-time domain $\mathbb{T}^2$ with two orthogonal coordinate variables. After connecting the pair of periodic boundaries of $\tau_2$, that is, $\mathbf{y}(\tau_1, 0) = \mathbf{y}(\tau_1, 2\pi)$, the domain $\mathbb{T}^2$ is transformed into a cylinder as shown in Fig. \ref{fig: tau_time}b. By further connecting the pair of periodic boundaries of $\tau_1$, that is, $\mathbf{y}(0, \tau_2) = \mathbf{y}(2\pi, \tau_2)$, a 2-torus yields as shown in Fig. \ref{fig: tau_time}c.

Changing the independent variable from $t$ to $\boldsymbol{\tau}$, system (2) is can be reformulated as $1^\text{st}$-order quasi-linear partial differential equations (PDEs)
\begin{equation}
	\sum_{i=1}^{d} \omega_i \frac{\partial \mathbf{y}(\boldsymbol{\tau})}{\partial \tau_i} = \mathbf{f}(\mathbf{y}(\boldsymbol{\tau}), \hat{\boldsymbol{\tau}}), \quad \mathbf{f}: \mathbb{R}^{2n \times 1} \times \mathbb{T}^e \mapsto \mathbb{R}^{2n \times 1},\label{EqPDEs}
\end{equation}
where $\hat{\boldsymbol{\tau}} = \boldsymbol{\Omega} \times t = [\tau_1, \dots, \tau_e]^\top$. The discretization techniques are usually adopted for the PDEs \eqref{EqA0} transformed from its $2^\text{nd}$ ODEs \eqref{2ODEs}, where displacements $\mathbf{q}(\boldsymbol{\tau}): \mathbb{T}^d \mapsto \mathbb{R}^{n \times 1}$ is usually approximated by a series of quasi-periodic shape functions with undetermined constant coefficients $\mathbf{Q}^d \in \mathbb{R}^{n\mathrm{U} \times 1}$. These quasi-periodic shape functions are essentially constructed by coupling $d$ sets of periodic functions, each associated with a hyper-time variable. Assume that the number of periodic functions with respect to $\tau_i$ is $\mathrm{U}_i$, there will be generally $\mathrm{U} = \prod_{i=1}^{d} \mathrm{U}_i$ quasi-periodic shape functions in discretization methods (see \ref{appA}). Determine these $n\mathrm{U}$ unknowns $\mathbf{Q}^d$ poses a significant challenge, especially for high-dimensional systems. 

\subsection{Multi-trajectory parametrization }\label{sec2.2}

\begin{figure}[ht] 
	\centering 
	\includegraphics[width=0.6\textwidth]{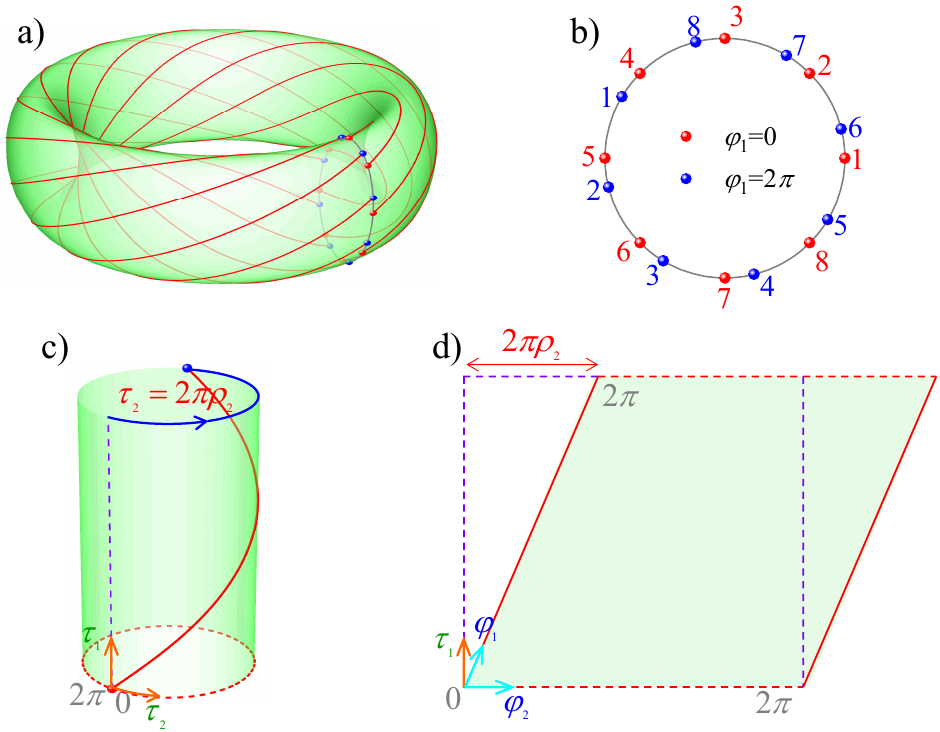} 
	\caption{\small The variant of hyper-time domain $\bar{\mathbb{T}}^d$ with $d = 2$: a) The multiple trajectories initialized at Poincaré section; b) The phase shift between initial and terminal points on trajectories; c) The representation of trajectories $\mathbf{z}(\varphi_1, 0) = \mathbf{y}(\varphi_1, \rho_2\varphi_1) $ and $\mathbf{z}(\varphi_1, 2\pi)= \mathbf{y}(\varphi_1,2\pi+\rho_2 \varphi_1)$ in hyper-time domain $\boldsymbol{\tau} \in \mathbb{T}^d$; d) The relationship between two hyper-time domains.} 
	\label{fig: phi_time_2}
\end{figure}

In this work, we define a new hyper-time domain $\boldsymbol{\varphi}=[\varphi_1,\varphi_2,\cdots,\varphi_d]^\top \in \mathbb{T}^d$ to parameterize $d$-torus as $\mathbf{z}(\varphi_1,\tilde{\boldsymbol{\varphi}})$ with $\tilde{\boldsymbol{\varphi}}=[\varphi_2,\varphi_3,\cdots,\varphi_d]^\top \in \mathbb{T}^{d-1}$. In this domain, $\mathbf{z}(0,\tilde{\boldsymbol{\varphi}})$ is initial first-order Poincaré section, and $\mathbf{z}(\varphi_1,\tilde{\boldsymbol{\varphi}}_s)$ represents a trajectory starting from initial point $\mathbf{z}(0,\tilde{\boldsymbol{\varphi}}_s)$.  For intuitive illustration, the case of $d = 2$ is exemplified in Fig. \ref{fig: phi_time_2}a, where the $2$-torus is represents as a collection of trajectories initialized on its Poincaré section (a limit cycle in Fig. \ref{fig: phi_time_2}b). By means of old domain $\boldsymbol{\tau}$, the trajectory $\mathbf{z}(\varphi_1,\tilde{\boldsymbol{\varphi}}_s)$ can be rewritten as
\begin{equation}
	\mathbf{z}(\varphi_1, \tilde{\boldsymbol{\varphi}}_s) := \mathbf{y}(\varphi_1, \tilde{\boldsymbol{\varphi}}_s + \boldsymbol{\rho} \varphi_1),\label{EqTraje}
\end{equation}

\begin{figure}[ht] 
 	\centering 
 	\includegraphics[width=0.75\textwidth]{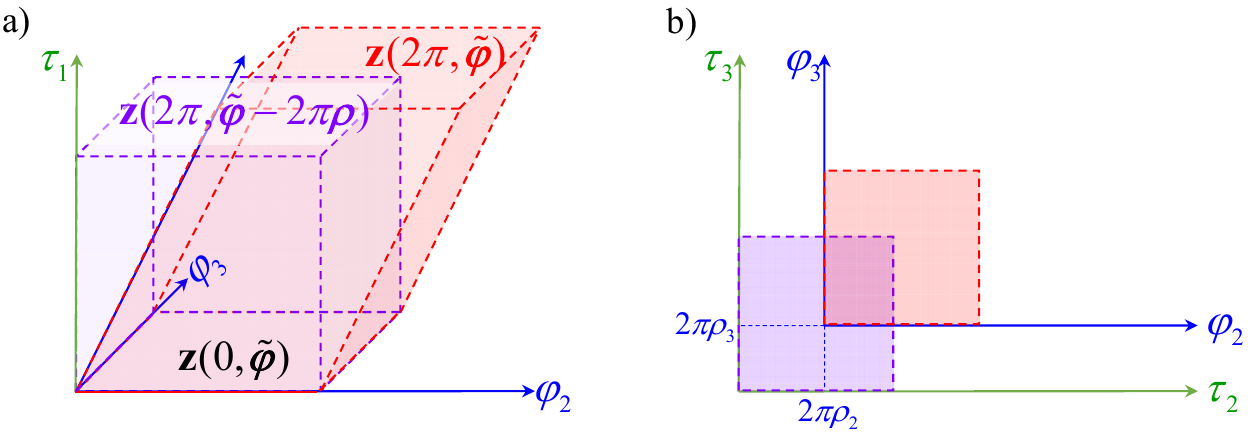} 
 	\caption{\small a): The relationship between the domains $\mathbb{T}^d$ and $\bar{\mathbb{T}}^d$ in the case with $d = 3$; b): The section on the two domains at $\tau_1 = \varphi_1 = 2\pi$.} 
 	\label{fig: phi_time_3}
 \end{figure}

Therefore, in Fig. \ref{fig: phi_time_2}c, for example, unfold the cylinder by cutting along the first trajectory $\mathbf{z}(\varphi_1, 0) := \mathbf{y}(\varphi_1, \rho_2\varphi_1)$,  the new domain of $(\varphi_1, \varphi_2) \in \bar{\mathbb{T}}^2$ yields as a parallelogram in Fig. \ref{fig: phi_time_2}d with respect to $(\tau_1, \tau_2) \in \mathbb{T}^2$. For the case with $d$ hyper-time variables, the domain $\boldsymbol{\tau} \in \mathbb{T}^d$ and $\boldsymbol{\varphi} \in \bar{\mathbb{T}}^d$ can be regarded respectively as a hyper-cubic and a hyper-parallelepiped, as shown in Fig. \ref{fig: phi_time_3}a for a case of $d = 3$. More importantly, following the definition of Eq. \eqref{EqTraje}, two types of boundary conditions can be obtained by $\mathbf{y}(\dots, 0, \dots) = \mathbf{y}(\dots, 2\pi, \dots)$ in domain $\boldsymbol{\tau} \in \mathbb{T}^d$
\begin{equation}
	\begin{aligned}
		\mathbf{z}(\varphi_1, \dots, 0, \dots) &= \mathbf{y}(\varphi_1, \dots, 0 + \rho_j \varphi_1, \dots) \\
		&= \mathbf{y}(\varphi_1, \dots, 2\pi + \rho_j \varphi_1, \dots) = \mathbf{z}(\varphi_1, \dots, 2\pi, \dots),\label{EqPerC}
	\end{aligned}
\end{equation}
and
\begin{equation}
	\mathbf{z}(0, \tilde{\boldsymbol{\varphi}}) = \mathbf{y}(0, \tilde{\boldsymbol{\varphi}}) = \mathbf{y}(2\pi, \tilde{\boldsymbol{\varphi}}) = \mathbf{z}(2\pi, \tilde{\boldsymbol{\varphi}} - 2\pi\boldsymbol{\rho}),\label{EqCoupC}
\end{equation}
where $\mathbf{z}(\varphi_1, \dots, 0, \dots)$ and $\mathbf{z}(\varphi_1, \dots, 2\pi, \dots)$ are a pair of periodic boundaries. For, example, coming back to the case with 2 hyper-time variables, $\mathbf{z}(\varphi_1,0)=\mathbf{y}(\varphi_1,\rho_2 \varphi_1)$ and $\mathbf{z}(\varphi_1, 2\pi)= \mathbf{y}(\varphi_1,2\pi+\rho_2 \varphi_1)$ are both red curve in Fig. \ref{fig: phi_time_2}c; Eq. \eqref{EqCoupC} reveal the coupling conditions between the initial state $\mathbf{z}(0, \tilde{\boldsymbol{\varphi}})$ and terminal state $\mathbf{z}(2\pi, \tilde{\boldsymbol{\varphi}})$ of the trajectories by using the phase shift $\tilde{\boldsymbol{\theta}} = 2\pi\boldsymbol{\rho}$ between them, as shown in Fig. \ref{fig: phi_time_2}d and Fig. 4b for the case of $d = 2$ and $d = 3$, respectively.

Based on the variable substitution and the chain rule, we obtain
\begin{equation}
	\begin{aligned}
		\omega_1 \frac{\partial \mathbf{y}(\boldsymbol{\tau})}{\partial \tau_1} &= \omega_1 \left( \frac{\partial \mathbf{z}(\varphi_1, \tilde{\boldsymbol{\varphi}})}{\partial \varphi_1} - \sum_{j=2}^{d} \rho_j \frac{\partial \mathbf{z}(\varphi_1, \tilde{\boldsymbol{\varphi}})}{\partial \varphi_j} \right) \\&= \omega_1 \frac{\partial \mathbf{z}(\varphi_1, \tilde{\boldsymbol{\varphi}})}{\partial \varphi_1} - \sum_{j=2}^{d} \omega_j \frac{\partial \mathbf{z}(\varphi_1, \tilde{\boldsymbol{\varphi}})}{\partial \varphi_j}; \\
		\omega_j \frac{\partial \mathbf{y}(\boldsymbol{\tau})}{\partial \tau_j} &= \omega_j \frac{\partial \mathbf{z}(\varphi_1, \tilde{\boldsymbol{\varphi}})}{\partial \varphi_j}, \quad j = 2, \dots, d.
	\end{aligned}\label{EqChainR}
\end{equation}
Thus the Eq. \eqref{EqPDEs} is transformed into
\begin{equation}
	\omega_1\mathbf{z}^{\prime}(\varphi_1,\tilde{\boldsymbol{\varphi}})=\mathbf{f}\left(\mathbf{z}(\varphi_1,\tilde{\boldsymbol{\varphi}}),\varphi_1,\hat{\boldsymbol{\varphi}}+\hat{\boldsymbol{\rho}}\varphi_1\right),\quad\mathbf{f}:\mathbb{R}^{2n\times1}\times\bar{\mathbb{T}}^e\mapsto\mathbb{R}^{2n\times1}\label{EqNewODEs}
\end{equation}
where $\hat{\boldsymbol{\varphi}} = [\varphi_2, \dots, \varphi_e]^\top$ and $\hat{{\boldsymbol\rho}} = [\rho_2, \dots, \rho_e]^\top$. $(\cdot)'$ denote derivative with respect to hyper-time variable $\varphi_1$. It is noticed that different from the complicated PDEs of Eq. \eqref{EqPDEs} in the hyper-time domain $\boldsymbol{\tau}$, the transformed system with the new parameterization variables is a set of ODEs with the same dimension as the original system \eqref{1ODEs}. By fixing $\tilde{\boldsymbol{\varphi}} = \tilde{\boldsymbol{\varphi}}_s$, any of trajectories $\mathbf{z}(\varphi_1, \tilde{\boldsymbol{\varphi}}_s)$ can be independently integrated from initial points $(0, \tilde{\boldsymbol{\varphi}}_s)$ to terminal points $(2\pi, \tilde{\boldsymbol{\varphi}}_s)$ based on the ODEs \eqref{EqNewODEs}. 

Accordingly, in this new hyper-time domain, the problem of determining the $d$-torus $\mathbf{z}(\varphi_1, \tilde{\boldsymbol{\varphi}})$ is simplified to searching for the $(d-1)$-torus $\mathbf{z}(0, \tilde{\boldsymbol{\varphi}})$ that satisfies the coupling conditions in Eq. \eqref{EqCoupC}. This is because the new parameterization method decouples $\varphi_1$ and $ \tilde{\boldsymbol{\varphi}}$. Since the problem is naturally reduced in dimension by one, the number of unknowns required to characterize the $(d-1)$-torus is also reduced. This lays the foundation for subsequent development of shooting method.

\section{Fourier series expansion-based shooting method}\label{sec3}

In this section, a Fourier series expansion-based shooting method, or \textit{FSE-Shooting} in short, is proposed to efficiently calculate the quasi-periodic solutions through parallel computing. First, a set of NAEs for the initial conditions is formulated by means of the coupling conditions in Sub-section \ref{sec3.1}. Then, Sub-section \ref{sec3.2} proposes the robust phase conditions to track the solution branches with unknown base frequencies through the one-parameter continuation. Next, in Sub-section \ref{sec3.3}, the implementation of FSE-Shooting is introduced by incorporating with the alternating frequency-time method (AFT) and time integration method (such as Newmark integration (NI)) . Finally, we also introduce the improvement of various methods in multi-trajectory parameterization based on FSE-Shooting. 

\subsection{Initial conditions and nonlinear algebraic equations}\label{sec3.1}

Given that the $2\pi$-periodicity of $(d-1)$-torus $\mathbf{z}(\varphi_1, \tilde{\boldsymbol{\varphi}})$ with respect to $\tilde{\boldsymbol{\varphi}}$ at any $\varphi_1$, as shown in Eq. \eqref{EqPerC}, $\mathbf{z}(\boldsymbol{\varphi})$ is approximated by the Fourier series expansion with variable-coefficients
\begin{equation}
	\mathbf{z}(\varphi_1, \tilde{\boldsymbol{\varphi}}) = (\mathbf{I}_{2n} \otimes \boldsymbol{\mathcal{H}}(\tilde{\boldsymbol{k}}, \tilde{\boldsymbol{\varphi}})) \mathbf{Z}(\tilde{\boldsymbol{k}}, \varphi_1),\label{EqFSE}
\end{equation}
where $\otimes$ represents Kronecker tensor product. $\mathbf{I}_{2n}$ is $2n \times 2n$ identity matrix. $\mathbf{Z}(\tilde{\boldsymbol{k}}, \varphi_1) \in \mathbb{R}^{2n\tilde{\mathrm{U}} \times 1}$ represents $2n\tilde{\mathrm{U}}$ variable-coefficients with respect to $\varphi_1$.  $\boldsymbol{\mathcal{H}}(\tilde{\boldsymbol{k}}, \tilde{\boldsymbol{\varphi}}) \in \mathbb{R}^{1 \times \tilde{\mathrm{U}}}$ is a series of quasi-periodic trigonometric functions, associated with $\tilde{\boldsymbol{\varphi}} = [\varphi_2, \dots, \varphi_d]^\top$
\begin{equation}
	\boldsymbol{\mathcal{H}}(\tilde{\boldsymbol{k}}, \tilde{\boldsymbol{\varphi}}) = [ 1, \dots, \cos(\tilde{\boldsymbol{k}}_l \tilde{\boldsymbol{\varphi}}), \sin(\tilde{\boldsymbol{k}}_l \tilde{\boldsymbol{\varphi}}), \dots ], \quad l = 1, \dots, \tilde{\mathrm{L}}^d,\label{EqQpTrig}
\end{equation}
where discretization parameter $\tilde{\boldsymbol{k}} = [k_2, \dots, k_d]$ represents the harmonic orders of Fourier series, with $\tilde{\boldsymbol{k}}_l$ being the $l$-th set of orders. The details for selecting $\tilde{\boldsymbol{k}}$ are introduced in \ref{appB}. Assume that $K_j = [0; \boldsymbol{k}_j] \in \mathbb{R}^{(\mathrm{L}_j + 1) \times 1}$ collects all elements of $|k_j|$ for $j = 2, \dots, d$, the value of $\tilde{\mathrm{U}}$ is actually equal to $\prod_{j=2}^{d} \mathrm{U}_j$, with $\mathrm{U}_j = 2\mathrm{L}_j + 1$. Hence, $\tilde{\mathrm{L}}^d$ in Eq. (11) is equal to $\tilde{\mathrm{U}} / 2 - 1$. 

Therefore, $\mathbf{z}(0, \tilde{\boldsymbol{\varphi}}) $ and $ \mathbf{z}(2\pi, \tilde{\boldsymbol{\varphi}} - 2\pi\boldsymbol{\rho})$ can be expressed via expansion \eqref{EqFSE} as
\begin{equation}
	\mathbf{z}(0, \tilde{\boldsymbol{\varphi}}) = (\mathbf{I}_{2n} \otimes \boldsymbol{\mathcal{H}}(\tilde{\boldsymbol{k}}, \tilde{\boldsymbol{\varphi}} )) \mathbf{Z}(\tilde{\boldsymbol{k}}, 0),\label{Eqz0}
\end{equation}
and
\begin{equation}
	\begin{aligned}
		\mathbf{z}(2\pi, \tilde{\boldsymbol{\varphi}} - 2\pi\boldsymbol{\rho})&= (\mathbf{I}_{2n} \otimes \boldsymbol{\mathcal{H}}(\tilde{\boldsymbol{k}}, \tilde{\boldsymbol{\varphi}} - 2\pi\boldsymbol{\rho}))\mathbf{Z}(\tilde{\boldsymbol{k}}, 2\pi)\\
		&= (\mathbf{I}_{2n} \otimes \boldsymbol{\mathcal{H}}(\tilde{\boldsymbol{k}}, \tilde{\boldsymbol{\varphi}})) (\mathbf{I}_{2n} \otimes \boldsymbol{\mathcal{R}}(\boldsymbol{\rho})) \mathbf{Z}(\tilde{\boldsymbol{k}}, 2\pi),
	\end{aligned}
	\label{Eqz2pi}
\end{equation}
with
\begin{equation}
	\boldsymbol{\mathcal{R}}(\boldsymbol{\rho}) = \mathrm{diag} \left( 1, \dots, \begin{bmatrix} \cos(2\pi\tilde{\boldsymbol{k}}_l \boldsymbol{\rho}) & -\sin(2\pi\tilde{\boldsymbol{k}}_l \boldsymbol{\rho}) \\ \sin(2\pi\tilde{\boldsymbol{k}}_l \boldsymbol{\rho}) & \cos(2\pi\tilde{\boldsymbol{k}}_l \boldsymbol{\rho}) \end{bmatrix}, \dots \right), \quad l = 1, \dots, \tilde{\mathrm{L}}^d
	,\label{EqRrho}
\end{equation}

Then, let us define a new set of Fourier coefficients $\hat{\mathbf{Z}}(\tilde{\boldsymbol{k}}, 2\pi)$ as
\begin{equation}
	\hat{\mathbf{Z}}(\tilde{\boldsymbol{k}}, 2\pi) = (\mathbf{I}_{2n} \otimes \boldsymbol{\mathcal{R}}(\boldsymbol{\rho})) \mathbf{Z}(\tilde{\boldsymbol{k}}, 2\pi)
	,\label{EqhatZ2pi}
\end{equation}
Here, $\mathbf{Z}(\tilde{\boldsymbol{k}}, 0)$ and $\hat{\mathbf{Z}}(\tilde{\boldsymbol{k}}, 2\pi)$ are Fourier coefficients of $\mathbf{z}(0, \tilde{\boldsymbol{\varphi}})$ and $\mathbf{z}(2\pi, \tilde{\boldsymbol{\varphi}} - 2\pi\boldsymbol{\rho})$ in the multi-frequency domain defined by $\boldsymbol{\mathcal{H}}(\tilde{\boldsymbol{k}}, \tilde{\boldsymbol{\varphi}})$, respectively. Thus, a set of NAEs of FSE-Shooting method is defined as ($\tilde{\boldsymbol{k}}$ is dropped for convenience)
\begin{equation}
	\begin{aligned}
		\mathbf{R}(\mathbf{Z}(0), \boldsymbol{\omega}) &= \hat{\mathbf{Z}}(2\pi) - \mathbf{Z}(0) \\
		&= (\mathbf{I}_{2n} \otimes \boldsymbol{\mathcal{R}}(\boldsymbol{\rho})) \mathbf{Z}(2\pi) - \mathbf{Z}(0).
	\end{aligned}
	\label{EqNAE_FSES}
\end{equation}
In order to define Fourier coefficients $\mathbf{Z}(0)$ as initial conditions of the FSE-Shooting, we design the following procedure to determine $\mathbf{Z}(2\pi)$
\begin{equation}
	\begin{aligned}
	   \mathbf{Z}(0) \overset{i\text{-DFT}^{d-1}}{\longrightarrow} \bar{\mathbf{z}}(0,\tilde{\boldsymbol{\varphi}}) \overset{\text{TI}}{\longrightarrow} \bar{\mathbf{z}}(2\pi,\tilde{\boldsymbol{\varphi}}) \overset{\text{DFT}^{d-1}}{\longrightarrow} \mathbf{Z}(2\pi)
		\end{aligned}
	\label{EqAFTNI}
\end{equation}
where $\bar{\mathbf{z}}(\varphi_1,\tilde{\boldsymbol{\varphi}})$ are $\tilde{\mathrm{S}} = \prod_{j=2}^{d} \mathrm{S}_j$ points selected in the sub-space of hyper-time domain. Here, $\mathrm{S}_j$ is the number of sampling points in the direction of $\varphi_j$. To avoid frequency aliasing associated with $\varphi_j$, the condition $\mathrm{S}_j \ge 2\mathrm{H}_j+1 \ge \mathrm{U}_j$ is required by Nyquist-Shannon sampling theorem, where $\mathrm{H}_j$ is the value of high order of $k_j$. This detail of procedure \eqref{EqAFTNI} will be introduced in Sub-section \ref{sec3.3}.

Noted that $\bar{\mathbf{z}}(0,\tilde{\boldsymbol{\varphi}})$ is initial condition in previous shooting methods \cite{R35,R76,R77}. However, these methods enforce $\tilde{\mathrm{S}}=\tilde{\mathrm{U}}$, a constraint that tends to induce frequency aliasing. However, the FSE-Shooting method not only reduces the number of unknowns to $2n\tilde{\mathrm{U}}\leq2n\tilde{\mathrm{S}}<n\mathrm{U}$, but also avoids frequency aliasing. In the special case that $d = 1$, the periodic solution $\mathbf{x}(t)$ can be transformed into $\mathbf{z}(\varphi_1) = \mathbf{Z}(\varphi_1)$. The Eq. \eqref{EqNAE_FSES} of FSE-Shooting degenerates into $\mathbf{R}(\mathbf{Z}(0), \omega_1) = \mathbf{Z}(2\pi) - \mathbf{Z}(0)$, which is actually the NAEs of the classical shooting method for periodic solution. 

\subsection{One-parameter continuation}\label{sec3.2} 

\subsubsection{Dimension deficit}\label{sec3.2.1} 

One-parameter continuation is usually used to compute roots of NAEs such as $\boldsymbol{\mathcal{R}}(\boldsymbol{\chi},p) = 0$, where $ \boldsymbol{\mathcal{R}}:\mathbb{R}^{N}\times\mathbb{R}\mapsto\mathbb{R}^{N} $ as the parameter $p$ evolves.  Here, $\boldsymbol{\chi}$ denotes the continuation variables and $p$ is continuation parameter. The two steps of one-parameter continuation reviewed in \ref{appC}, namely tangent prediction and orthogonal corrections, are usually used to ensure the continuation process.

From the analysis of FSE-Shooting or \textit{m}-VCF and  \textit{s}-CCF in \ref{appA}, all of them construct a set of NAEs with form $\mathbf{R}(\mathbf{C},\boldsymbol{\omega},\lambda) = \mathbf{0}$, where $\mathbf{R}: \mathbb{R}^{N_C}\times\mathbb{R}^d\times\mathbb{R}^a \mapsto \mathbb{R}^{N_C}$. For example, Eq. \eqref{EqNAE_FSES} in FSE-Shooting and Eq. \eqref{EqA2} in Appendix. Here, $\lambda$ represents the system parameters excluding the $e$ excitation fundamental frequencies $\mathbf{\Omega}$. Noted that there are $d-e$ base frequencies $\boldsymbol{\omega}_u$ not defined by the system or excitation. The $d-e$ phase conditions $PC_i(\mathbf{C}) = 0, \, i=e+1,e+2,\cdots,d$ are usually added to fix the phase of $d$-torus with respect to unknown base frequencies, which will be detailed in next sub-section. Therefore, we usually construct the following continuation problem to track $d$-torus
\begin{equation}
	\boldsymbol{\mathcal{R}}([\mathbf{C};\boldsymbol{\omega}_u],p) = \begin{cases}
		\mathbf{R}(\mathbf{C},\boldsymbol{\omega},\lambda) = \mathbf{0} \\
		PC_{i}(\mathbf{C}) = 0, \quad i=e+1,\cdots,d\\		
	\end{cases}\label{EqCont1}
\end{equation}
where $[\mathbf{C};\boldsymbol{\omega}_u]$ collects $N_C+d-e$ continuation variables; the continuation parameter $p$ is selected from $\mathbf{\Omega}$ and $\lambda$, and those not selected must be fixed. In addition, when one of $\mathbf{\Omega}$ is chosen as $p$, for example $\omega_1$, a more robust continuation problem is defined as
\begin{equation}
	\boldsymbol{\mathcal{R}}([\mathbf{C};\boldsymbol{\omega}_u;\omega_2;\cdots;\omega_e],\omega_1) = \begin{cases}
		\mathbf{R}(\mathbf{C},\boldsymbol{\omega}) = \mathbf{0} \\
		PC_{i}(\mathbf{C}) = 0, \quad i=e+1,\cdots,d\\
		FC_{j}(\omega_j,\omega_1) = 0, \quad i=2,\cdots,e\\
	\end{cases}\label{EqCont2}
\end{equation}
which has $N_C+d-1$ continuation variables. So, $e-1$ frequency conditions, for example $FC_j = \omega_j - \rho_j\omega_1 = 0$, are added to ensure the continuation is well-posed.

\subsubsection{Phase conditions}\label{sec3.2.2} 

Given the first-order ODEs \eqref{1ODEs} in torus coordinates, they can be expressed as $\dot{\mathbf{x}}(\boldsymbol{\theta})=\mathbf{f}(\mathbf{x},\hat{\boldsymbol{\theta}})$, where $\hat{\boldsymbol{\theta}}=[\theta_1,\cdots,\theta_e]^\top$. Suppose the $d$-torus $\mathbf{x}(\boldsymbol{\theta})$ satisfies this system, the $d$-torus $\mathbf{x}(\boldsymbol{\theta}+\boldsymbol{\theta}_a)$ with phase shift $\boldsymbol{\theta}_a$ must be a solution of $\dot{\mathbf{x}}(\boldsymbol{\theta}+\boldsymbol{\theta}_a)=\mathbf{f}(\mathbf{x},\hat{\boldsymbol{\theta}}+\hat{\boldsymbol{\theta}}_a)$. However, setting $\hat{\boldsymbol{\theta}}_a = \mathbf{0}$, i.e., $\boldsymbol{\theta}_a=[\mathbf{0};\theta_{e+1,a};\cdots;\theta_{d,a}]$, $\mathbf{x}(\boldsymbol{\theta}+\boldsymbol{\theta}_a)$ can also be regarded as a solution to the original system. That is to say, when $d>e$, the phases of $d$-torus with respect to unknown base frequencies are undetermined. Thus the coefficients $\mathbf{C}$ are functions with respect to $\theta_{i,a},\,i=e+1,\cdots,d$.

The phase conditions are introduced to fix these phases, or equivalently, to specify the characteristics of $d$-torus at zero phases, i.e., $\theta_{i,a}:=\theta_{i,0}=0$. The most common variants are derivative condition \cite{R61}, Poincaré condition \cite{R74} and integral Poincaré condition \cite{R35}
\begin{equation}
	PC_i =  \left.  \frac{\partial x_m}{\partial \theta_i} \right|_{\boldsymbol{\theta}_0} = 0,
	\label{EqPC1}
\end{equation}
\begin{equation}
	PC_i =  \left.(\mathbf{x}-\mathbf{x}^*)^\top \frac{\partial \mathbf{x}^*}{\partial \theta_i}\right|_{\boldsymbol{\theta}_0} = 0,
	\label{EqPC2}
\end{equation}
and
\begin{equation}
	PC_i = \int_{0}^{2\pi} \mathbf{x}^\top   \frac{\partial \mathbf{x}^*}{\partial \theta_i}  \mathrm{d}\boldsymbol{\theta} = 0,
	\label{EqPC3}
\end{equation}
where $x_m$ is the $m$-th variable of $\mathbf{x}$; $\mathbf{x}^*$ is a known solution in the neighborhood of $\mathbf{x}$, i.e., the previous solution on continuation curve. Among these phase conditions, the convergence of derivative condition is affected by $m$ and the initial guess $\mathbf{C}_{j+1}^0$; the Poincaré condition is also influenced by $\mathbf{C}_{j+1}^0$, and the integral Poincaré condition is the most robust and it eventually takes the form
\begin{equation}
	PC_i = \mathbf{C}^\top \mathbf{O}\mathbf{C}^* = 0,
\end{equation}
where $\mathbf{O}\in \mathbb{R}^{N_C\times N_C}$ is transformed matrix, and the two vectors satisfy the orthogonality relation $\mathbf{C} \perp \mathbf{O}\mathbf{C}^*$. However, it suffers from two disadvantages: $\mathbf{C}^*$ needs to be stored extra, and no reference $\mathbf{C}^*$ is available for the first solution during continuation.

To overcome the aforementioned drawbacks, a robust phase condition was proposed in our previous work \cite{R41,R66}, which is also extended for FSE-Shooting. In tangent prediction and orthogonal corrections, we also find two orthogonal relations associated with the conventional phase condition
\begin{equation}
		\begin{aligned}
			 \partial_{\mathbf{C}} PC_i^\top |_{(\boldsymbol{\chi}_j, p_j)} \perp  \Delta \mathbf{C}_j  ;\quad 
			\partial_{\mathbf{C}} PC_i^\top|_{(\boldsymbol{\chi}_{j+1}^k, p_{j+1}^k)} \perp \delta \mathbf{C}_j^{k+1}
		\end{aligned}
		\label{Eq2Op}
\end{equation}
where the second holds when $\|PC_i\|_2 < \varepsilon$. Therefore, If we directly define a vector related to $\mathbf{C}$ that is always orthogonal to both $\Delta \mathbf{C}_j$ and $\delta \mathbf{C}_j^{k+1}$ during continuation, there is no need to store $\mathbf{C}^*$ additionally. Inspired by integral Poincaré condition, we propose the following phase condition for FSE-Shooting method
\begin{equation}
	PC_i = \int_{0}^{2\pi} \cdots \int_{0}^{2\pi} \left[ \frac{\partial \mathbf{z}(\varphi_1, \tilde{\boldsymbol{\varphi}})^\top}{\partial \varphi_i} \mathbf{z}^*(\varphi_1, \tilde{\boldsymbol{\varphi}}) \right]_{\varphi_1 = 0} \mathrm{d}\varphi_2 \cdots \mathrm{d}\varphi_d \equiv 0,
	\label{EqOPC}
\end{equation}
By using Fourier-Galerkin procedure in $\tilde{\boldsymbol{\varphi}} \in \tilde{\mathbb{T}}^{d-1}$ 
\begin{equation}
	PC_i = \begin{cases}
		\mathbf{Z}'(0) \perp \mathbf{Z}^*(0)  & i = 1; \\
		\left[ (\mathbf{I}_{2n} \otimes \boldsymbol{\nabla}_{\varphi_i}) \mathbf{Z}(0) \right] \perp \mathbf{Z}^*(0) & i \neq 1,
	\end{cases}
	\label{Eq18}
\end{equation}
where $\mathbf{Z}^*(0):=\Delta \mathbf{Z}(0)_j$ or $\mathbf{Z}^*(0):=\delta \mathbf{Z}(0)_j^{k+1}$ are defined as the Fourier coefficients of $\mathbf{z}^*(\varphi_1, \tilde{\boldsymbol{\varphi}})$ in Eq. \eqref{EqOPC}. $\mathbf{Z}'(0)$ represents the corresponding Fourier coefficients of $\mathbf{z}'(\varphi_1, \tilde{\boldsymbol{\varphi}})$ with $\varphi_1 = 0$. And the constant matrix $\boldsymbol{\nabla}_{\varphi_i} \in \mathbb{R}^{\tilde{\mathrm{U}} \times \tilde{\mathrm{U}}}$ is governed by
\begin{equation}
	\boldsymbol{\nabla}_{\varphi_i} = \mathrm{diag} \left(0 , \dots, \begin{bmatrix}
		0 & k_{i,l}  \\
		-k_{i,l} & 0  \\
	\end{bmatrix},\cdots \right), \quad i = 2, \dots, d, \quad l = 1, \dots, \tilde{\mathrm{L}}^d
	\label{Eq19}
\end{equation}

It is worth emphasizing that the proposed phase condition is not a nonlinear equation for $\mathbf{Z}(0)$, but an orthogonal geometric condition that locks the zero phase by enforcing orthogonality between the two vectors. By eliminating reliance on state variables, initial values, and the previous solution, the robustness of proposed phase condition has been demonstrated in the discretization method \cite{R41,R66}.

\subsection{Implementation of FSE-Shooting}\label{sec3.3} 

\begin{figure}[ht] 
	\centering 
	\includegraphics[width=0.85\textwidth]{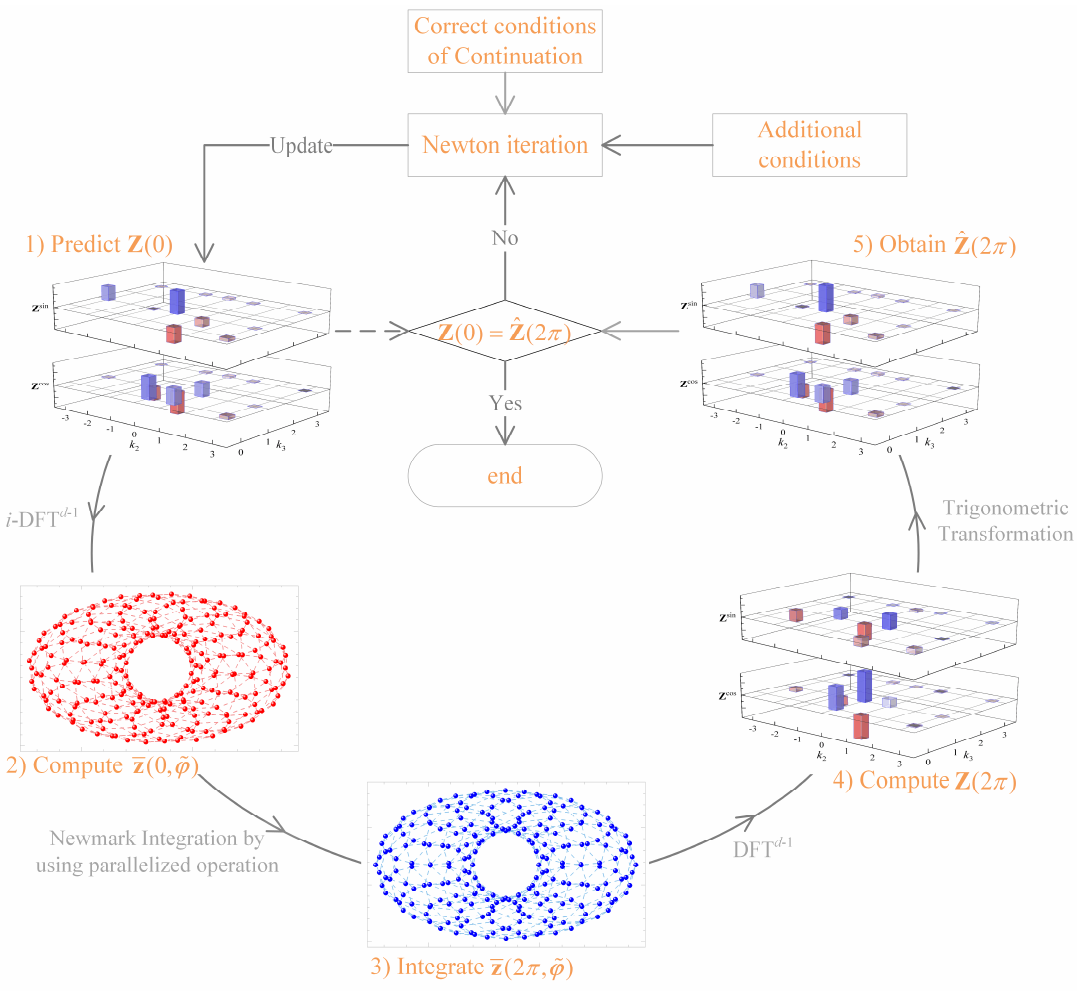} 
	\caption{\small The flowchart of the Fourier series expansion-based shooting method (an example for quasi-periodic solution with 3 base frequencies).} 
	\label{fig: FSE-Shooting}
\end{figure}

Based on Eqs. \eqref{EqNAE_FSES} and \eqref{EqAFTNI}, the implementation of the FSE-Shooting method is introduced in Fig. \ref{fig: FSE-Shooting}, which contains the following 6 steps:

1): Initialize or update the initial Fourier coefficients $\mathbf{Z}(0)$ of $(d-1)$-torus $\mathbf{z}(0, \tilde{\boldsymbol{\varphi}})$ (Sometime, also predict unknown base frequencies $\boldsymbol{\omega}_u$);

2): By using the $d-1$ dimensional inverse discrete Fourier transforms ($i$-DFT$^{d-1}$) of the alternating frequency-time method (AFT), compute sampling points $\bar{\mathbf{z}}(0, \tilde{\boldsymbol{\varphi}})$, as governed by
\begin{equation}
	\bar{\mathbf{z}}(0, \tilde{\boldsymbol{\varphi}}) = (\mathbf{I}_{2n} \otimes \boldsymbol{\Gamma}) \mathbf{Z}(0)
	\label{Eq20}
\end{equation}
where $\boldsymbol{\Gamma} \in \mathbb{R}^{\tilde{\mathrm{S}} \times \tilde{\mathrm{U}}}$ is the constant matrix of $i$-DFT$^{d-1}$, introduced in \ref{appD};

3): By dividing $\varphi_1 \in [0, 2\pi)$ into $\mathrm{S}_1$ sub-intervals, the time integration, such as Newmark integration, is used to parallelly compute the terminal points $\bar{\mathbf{z}}(2\pi, \tilde{\boldsymbol{\varphi}}) $ of trajectories started at $\bar{\mathbf{z}}(0, \tilde{\boldsymbol{\varphi}}) $. Also, the derivatives of $\partial_{\bar{\mathbf{z}}(0, \tilde{\boldsymbol{\varphi}})} \bar{\mathbf{z}}(2\pi, \tilde{\boldsymbol{\varphi}}) \in \mathbb{R}^{2n\tilde{\mathrm{S}} \times 2n\tilde{\mathrm{S}}}$ and $\partial_{\omega_1} \bar{\mathbf{z}}(2\pi, \tilde{\boldsymbol{\varphi}}) \in \mathbb{R}^{2n\tilde{\mathrm{S}} \times 1}$ are also obtained simultaneously, as detailed in \ref{appE}. More importantly, like the shooting method for periodic solution, the by-product $\partial_{\bar{\mathbf{z}}(0, \tilde{\boldsymbol{\varphi}})} \bar{\mathbf{z}}(2\pi, \tilde{\boldsymbol{\varphi}}) $ of FSE-Shooting method can be also directly used to compute the Lyapunov exponents to assess the stabilities of quasi-periodic solutions, as shown in \ref{appF};

4): By using the DFT$^{d-1}$ of AFT, the Fourier coefficients $\mathbf{Z}(2\pi)$ are transformed from $\bar{\mathbf{z}}(2\pi, \tilde{\boldsymbol{\varphi}})$:
\begin{equation}
	\mathbf{Z}(2\pi) = (\mathbf{I}_{2n} \otimes \boldsymbol{\Gamma}^{-1}) \bar{\mathbf{z}}(2\pi, \tilde{\boldsymbol{\varphi}})
	\label{Eq21}
\end{equation}
where $\boldsymbol{\Gamma}^{-1} \in \mathbb{R}^{\tilde{\mathrm{U}} \times \tilde{\mathrm{S}}}$ is the constant matrix of DFT$^{d-1}$ and the detail is also introduced in \ref{appD};

5): By using Trigonometric transformation $\boldsymbol{\mathcal{R}}(\boldsymbol{\rho})$, the Fourier coefficients $\hat{\mathbf{Z}}(2\pi)$ of $\mathbf{z}(2\pi, \tilde{\boldsymbol{\varphi}} - 2\pi\boldsymbol{\rho})$ are obtained by $\hat{\mathbf{Z}}(2\pi) = (\mathbf{I}_{2n} \otimes \boldsymbol{\mathcal{R}}(\boldsymbol{\rho})) \mathbf{Z}(2\pi)$;

6): Check whether $\mathbf{Z}(0)$ equals $\hat{\mathbf{Z}}(2\pi)$. If so, the initial conditions of the quasi-periodic solution are determined. Conversely, use the Newton-Raphson method to update initial conditions $\mathbf{Z}(0)$ and $\boldsymbol{\omega}$, and return back to step 1.

And the Jacobian associated with the NAEs \eqref{EqNAE_FSES} are computed by:
\begin{equation}
	\partial_{\mathbf{Z}(0)} \mathbf{R} = (\mathbf{I}_{2n} \otimes \boldsymbol{\mathcal{R}}(\boldsymbol{\rho}) \boldsymbol{\Gamma}^{-1}) \frac{\partial \bar{\mathbf{z}}(2\pi, \tilde{\boldsymbol{\varphi}})}{\partial \bar{\mathbf{z}}(0, \tilde{\boldsymbol{\varphi}})} (\mathbf{I}_{2n} \otimes \boldsymbol{\Gamma}) - \mathbf{I}_{2n\tilde{\mathrm{U}}},
	\label{Eq22}
\end{equation}
\begin{equation}
	\partial_{\omega_i} \mathbf{R} = \begin{cases}
		(\mathbf{I}_{2n} \otimes \boldsymbol{\mathcal{R}}(\boldsymbol{\rho}) \boldsymbol{\Gamma}^{-1}) \frac{\partial \bar{\mathbf{z}}(2\pi, \tilde{\boldsymbol{\varphi}})}{\partial \omega_i} + (\mathbf{I}_{2n} \otimes \partial_{\omega_i} \boldsymbol{\mathcal{R}}(\boldsymbol{\rho})) \mathbf{Z}(2\pi), & i = 1,2,\cdots,e; \\
		(\mathbf{I}_{2n} \otimes \partial_{\omega_i} \boldsymbol{\mathcal{R}}(\boldsymbol{\rho})) \mathbf{Z}(2\pi), & i = e+1,\cdots,d,
	\end{cases}
	\label{Eq23}
\end{equation}
where
\begin{equation}
	\frac{\partial \boldsymbol{\mathcal{R}}(\boldsymbol{\rho})}{\partial \omega_i} = \frac{2\pi}{\omega_i}
	\mathrm{diag} \left(0 , \dots, \boldsymbol{\Upsilon}_{i,l} ,\cdots \right)\boldsymbol{\mathcal{R}}(\boldsymbol{\rho})
	\label{Eq24}
\end{equation}
with $\boldsymbol{\Upsilon}_{1,l} = -\begin{bmatrix} 0 & -\tilde{\boldsymbol{k}}_l \boldsymbol{\rho} \\ \tilde{\boldsymbol{k}}_l \boldsymbol{\rho} & 0 \end{bmatrix}$ and $\boldsymbol{\Upsilon}_{i,l} = \begin{bmatrix} 0 & -k_{i,l} \\ k_{i,l} & 0 \end{bmatrix}$, $i = 2, \dots, d,\, l = 1, \dots, \tilde{\mathrm{L}}^d$. Moreover, the partial-derivative $\partial_p \mathbf{R}$ can be assessed by the finite difference method when $p \neq \omega_i$.

\subsection{Improvement of other methods in multi-trajectory parameterization}\label{sec3.4} 

Shooting methods proposed in \cite{R35,R76,R77} focus on searching initial points $\bar{\mathbf{z}}(0,\tilde{\boldsymbol{\varphi}})$, FSE-Shooting coupled with AFT can be rewritten as
\begin{equation}
	\begin{aligned}
		\mathbf{R}(\bar{\mathbf{z}}(0,\tilde{\boldsymbol{\varphi}}), \boldsymbol{\omega}) = (\mathbf{I}_{2n} \otimes \boldsymbol{\Gamma} \boldsymbol{\mathcal{R}}(\boldsymbol{\rho}) \boldsymbol{\Gamma}^{-1}) \bar{\mathbf{z}}(2\pi,\tilde{\boldsymbol{\varphi}}) - \bar{\mathbf{z}}(0,\tilde{\boldsymbol{\varphi}})=\mathbf{0},
	\end{aligned}
	\label{EqFSES_T}
\end{equation}
where the matrices $\boldsymbol{\Gamma}$ and $\boldsymbol{\Gamma}^{-1}$ do not require $\tilde{\mathrm{S}}=\tilde{\mathrm{U}}$,
thus avoiding frequency aliasing. Eq. \eqref{EqFSES_T} can be referred to as FSE-Shooting in the hyper-time domain, and the proposed method is FSE-Shooting in the multi-frequency domain. 

These trajectories can also be determined by discretization method, such as collocation method \cite{R22,R50}, where undetermined coefficients are defined as $\mathbf{C}:=\bar{\mathbf{z}} = [\bar{\mathbf{z}}(\varphi_{1,0},\tilde{\boldsymbol{\varphi}});\bar{\mathbf{z}}(\varphi_{1,1},\tilde{\boldsymbol{\varphi}});\cdots;\bar{\mathbf{z}}(\varphi_{1,\mathrm{S}_1},\tilde{\boldsymbol{\varphi}})]$, with $\varphi_{1,k} = k \Delta \varphi_1 = 2k\pi/\mathrm{S}_1$. The number of $\bar{\mathbf{z}}$ is $ 2n(\mathrm{S}_1+1)\tilde{\mathrm{S}}$ and the NAEs consist of two parts
\begin{equation}
	\mathbf{R}(\bar{\mathbf{z}}, \boldsymbol{\omega}) = \begin{cases}
		\mathbf{R}(\bar{\mathbf{z}}(0,\tilde{\boldsymbol{\varphi}}), \boldsymbol{\omega}) =\mathbf{0}\\
	\mathbf{R}(\bar{\mathbf{z}}(\bar{\boldsymbol{\varphi}}_1,\tilde{\boldsymbol{\varphi}}_s), \boldsymbol{\omega}) =\mathbf{0}, \quad  s = 1,2,\cdots,\tilde{\mathrm{S}}	
	\end{cases}
	\label{EqFSEDM}
\end{equation}	
where $\mathbf{R}(\bar{\mathbf{z}}(\bar{\boldsymbol{\varphi}}_1,\tilde{\boldsymbol{\varphi}}_s), \boldsymbol{\omega}) =\mathbf{0}$ has $2n\mathrm{S}_1$ equations.

It is not difficult to observe that FSE-Shooting has the fewest unknowns among all methods, and with the aid of the AFT method, it avoids frequency aliasing. However, the efficiency of numerical method is not only influenced by the dimension of NAEs; the construction of Jacobian matrix also plays an equally crucial role. The parallelization operation in FSE-Shooting helps improve this construction speed. In the next section, the efficiency of proposed method will be verified. For comparison, this paper adopts the VCHB constructed by \textit{m}-VCF, which is the most efficient approach in multi-period parameterization.

\section{Application to nonlinear finite element systems}\label{sec4}

In this section, all simulations have been performed on a computer with 32-Core Processor, whose base speed is 3.10 GHz. The programming language is MATLAB and the version R2022a is used.

All three finite element systems belong to the Eq. \eqref{EqCont2} introduced in Sub-section \ref{sec3.2.1}, where $\omega_1$ is selected as the continuation parameter. Three cases with different $e$ and $d$ are discussed, i.e., $d = e = 2$ in Sub-section \ref{sec4.1.1} with 237 DOFs and \ref{sec4.3} with 1872 DOFs; $d = e = 3$ in Sub-section \ref{sec4.1.2} with 33 DOFs; and $d = 2$, $e = 1$ in Sub-section \ref{sec4.2} with 94 DOFs.

\subsection{A quasi-periodically-forced von Kármán beam}\label{sec4.1} 

\begin{figure}[ht] 
	\centering 
	\includegraphics[width=0.75\textwidth]{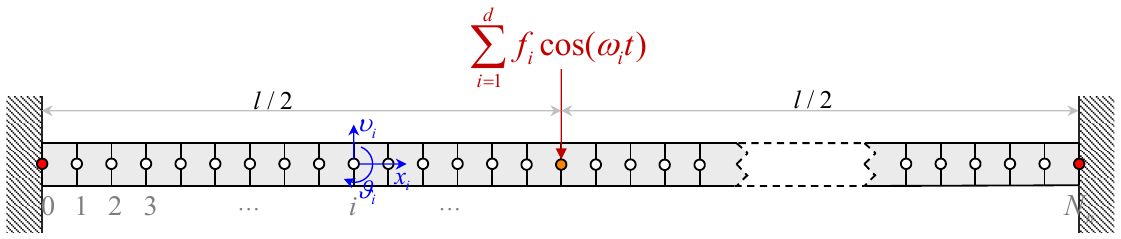} 
	\caption{\small A Finite-element clamped-clamped von Kármán beam with quasi-periodic excitation at its midspan.} 
	\label{fig: CC-VKB}
\end{figure}

The considered test case is that of a clamped-clamped von Kármán beam with quasi-periodic excitation at its midspan, as depicted in Fig. \ref{fig: CC-VKB}, where the axial displacement $x_i$, the transverse displacement $v_i$ and the rotation angle $\vartheta_i$ are introduced at each node by using classic finite-element method. Assume that the beam is divided into $N_n$ elements, there are $3N_n - 3$ degrees of freedom after boundary conditions are imposed. $\sum_{i=1}^{d} f_i \cos(\omega_i t)$ represents the quasi-periodic excitations at the midspan of beam, and the base frequencies have the ratios of $\rho_2 = \omega_2 / \omega_1 = 1 / \sqrt{2}$ and $\rho_3 = \omega_3 / \omega_1 = 1 / \sqrt{5}$. The parameters of the von Kármán beam system are summarized in Table \ref{tab1}.

\begin{table}[htbp]
	\centering
	\small
	\caption{ Physical parameters of the clamped-clamped von Kármán beam system}
	\label{tab1}
		\begin{tabular}{ll}
			\toprule
			Physical parameter & Value \\
			\midrule
			Young's modulus $E$ [kPa] & $4.5 \times 10^7$ \\
			Density $\rho$ [kg/mm$^3$] & $1780 \times 10^{-9}$ \\
			Height of beam $h$ [mm] & $10$ \\
			Width of beam $b$ [mm] & $10$ \\
			Length of beam $l$ [mm] & $2700$ \\
			Poisson's ratio $\nu$ & $0.3$ \\
			Mass-proportional damping coefficient $\alpha$ [s] & $0$ \\
			Stiffness-proportional damping coefficient $\beta$ [s$^{-1}$] & $2 / 9 \times 10^{-4}$ \\
			Amplitude of excitation $f_1, f_2, f_3$ [N] & $1, 1, 1$ \\
			\bottomrule
		\end{tabular}
\end{table}

The clamped-clamped von Kármán beam can be modeled by YetAnotherFEcode \cite{R67}, integrated into SSMTool \cite{R68}:
\begin{equation}
	\mathbf{M} \ddot{\mathbf{q}} + \mathbf{D} \dot{\mathbf{q}} + \mathbf{K} \mathbf{q} + \mathbf{f}_{nl}(\mathbf{q}) = \sum_{i=1}^{d} \boldsymbol{f}_i \cos(\omega_i t)
	\label{Eq25}
\end{equation}
where $\mathbf{q} \in \mathbb{R}^{(3N_n - 3) \times 1}$ is the vector of generalized displacements in form of $\mathbf{x} = [x_1; v_1; \vartheta_1; \dots; x_{N_n - 1}; v_{N_n - 1}; \vartheta_{N_n - 1}]$. $\mathbf{M}$, $\mathbf{D}$, $\mathbf{K}$ are the mass, the damping and the stiffness matrices. $\mathbf{f}_{nl}(\mathbf{q})$ is the vector of the nonlinear factors introduced by the axial stretching due to large transverse deflections. And the damping matrix is set as $\mathbf{D} = \alpha \mathbf{M} + \beta \mathbf{K}$, where $\alpha$, $\beta$ represent the mass-proportional damping coefficient and stiffness-proportional damping coefficient, respectively. $\boldsymbol{f}_i$ represents the vector of the amplitude of excitation with $\boldsymbol{f}_{i, 3N_n / 2 - 1} = f_i$ and zero otherwise. Here, $f_i$ is the amplitude of excitation forced at midspan of beam.

\subsubsection{Quasi-periodic forcing with 2 base frequencies}\label{sec4.1.1}

The case under the quasi-periodic excitation with $e = 2$ and quasi-periodic solutions with $d = 2$ is firstly focused, the number of elements of beam is selected as 80. There are 237 DOFs in system \eqref{Eq25}, also 474 state variables in $1^\text{st}$ ODEs. In particular, the equations of motion \eqref{Eq25} are rewritten in the first-order form \eqref{1ODEs} with both displacement and velocity treated as state variables. $\omega_1$ is selected as the continuation parameter.

\begin{figure}[ht] 
	\centering 
	\includegraphics[width=0.8\textwidth]{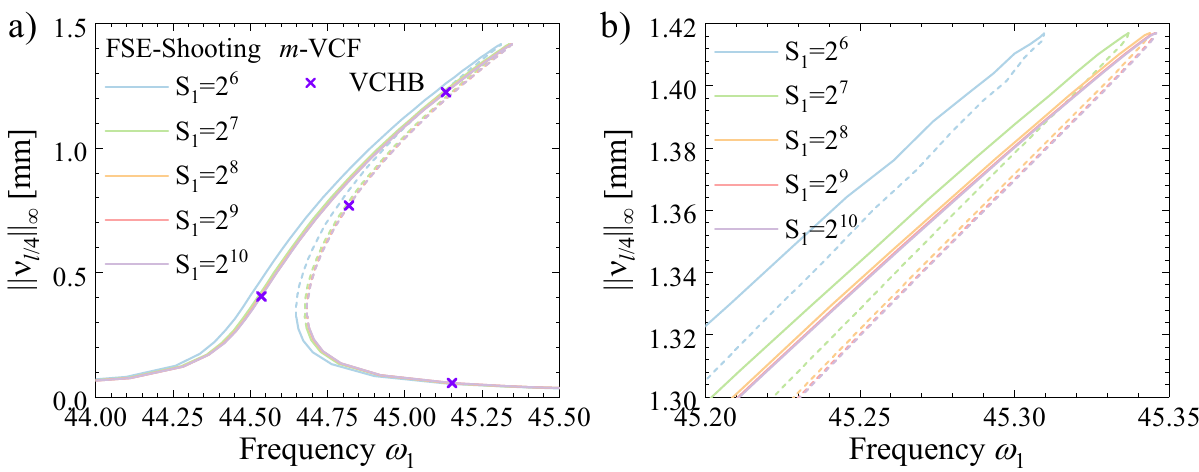} 
	\caption{\small a): The FRCs in the amplitude of transverse displacement at $l/4$ for clamped-clamped von Kármán beam under the quasi-periodic excitation with two base frequencies by using the proposed method (FSE-Shooting) and previous method (\textit{m}-VCF); b): Zoom of peaks.} 
	\label{fig: FRCs d=e=2}
\end{figure}

To investigate the influence of the step $\Delta \varphi_1$ in the Newmark integration on the results, the frequency response curves (FRCs) in the amplitude of transverse displacement at $l/4$ for clamped-clamped von Kármán beam with $\mathrm{S}_1 = 2^6$, $2^7$, $2^8$, $2^9$, $2^{10}$ are shown in Fig. \ref{fig: FRCs d=e=2}. The other parameters of FSE-Shooting method are defined as $K_2 = [0, 1, \dots, 4, 5]^\top \in \mathbb{R}^{6 \times 1}$, $\mathrm{S}_2 = 2^5$ and $\mathrm{\tilde{U}} = 11$. As a consequence, there are 5214 unknown coefficients $\mathbf{Z}(0)$ in FSE-Shooting \eqref{EqNAE_FSES}. The frequency condition $g_2 = \omega_2 - \rho_2 \omega_1 = 0$ is added into the continuation. In Fig. \ref{fig: FRCs d=e=2}, the slide curves represent the stable quasi-periodic solutions, and dotted curves denote unstable ones. As the parameter $\mathrm{S}_1$ increases, the FRCs gradually converge. 

\begin{figure}[ht] 
	\centering 
	\includegraphics[width=0.8\textwidth]{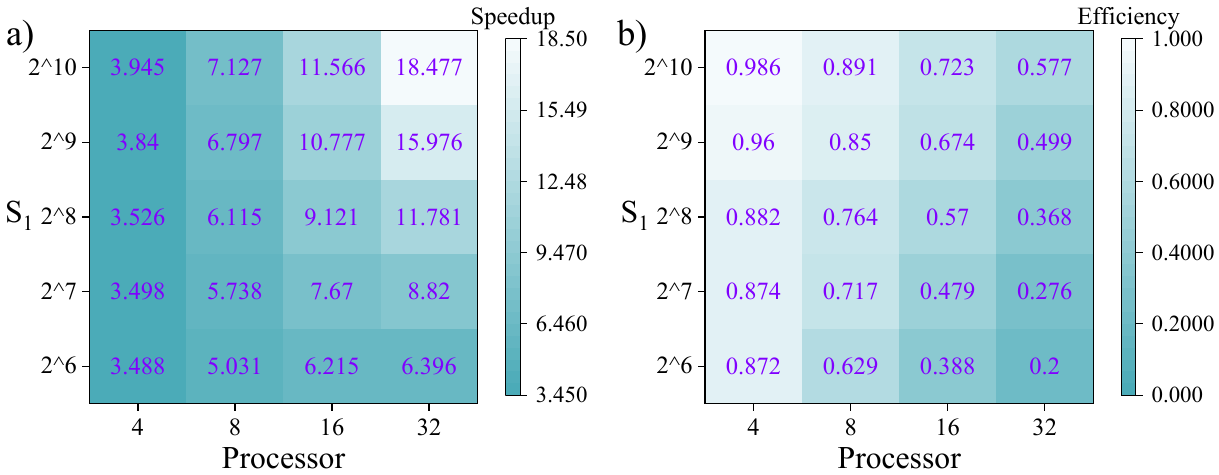} 
	\caption{\small Investigation of computation cost per iteration by using parallel with respect to step parameter $\mathrm{S}_1$ and Processors of CPU in Newmark integration, a): The speedup; b): The efficiency.} 
	\label{fig: Speedup & efficiency}
\end{figure}

Then, by using parallel in Newmark integration, the computation costs per iteration of the continuation are concluded in Table \ref{tab2}. As the number $\mathrm{S}_1$ of integration increases, the computation cost shows an increasing trend, without parallel. By introducing the parallel operations, the computational cost decreases significantly as the number of Processors of CPU increases. The speedup and efficiency of parallel \cite{R69} are shown in Fig. \ref{fig: Speedup & efficiency}, where the former is the ratio of serial execution time to parallel execution time and the latter is defined as the ratio of speedup to the number of Processors. The results show that the highest value of speedup 18.477 occurs at $\mathrm{S}_1 = 2^{10}$ and 32 Processors, the highest value of efficiency 0.986 occurs at $\mathrm{S}_1 = 2^{10}$ and 4 Processors. Although the increase in the number of Processors of CPU leads to a certain decline in parallel efficiency, it still contributes to reducing the computational cost, i.e., the speedup ratio is increasing.

\begin{table}[htbp]
	\centering
	\caption{\small The computation costs per iteration [s] of continuation by using parallel with respect to step parameter $\mathrm{S}_1$ and Processors of CPU in Newmark integration}
	\label{tab2}
	\begin{tabular}{cccccc}
		\toprule
		\multirow {2}{*}{Number of $\mathrm{S}_1$} & \multicolumn{5}{c}{Processors of CPU} \\
		\cmidrule(lr){2-6}
		& w/o parallel & 4 & 8 & 16 & 32 \\
		\midrule
		$2^6$ & 49.307 & 14.133 & 9.801 & 7.934 & 7.709 \\
		$2^7$ & 86.166 & 24.633 & 15.015 & 11.235 & 9.769 \\
		$2^8$ & 160.232 & 45.440 & 26.201 & 17.567 & 13.601 \\
		$2^9$ & 331.246 & 86.267 & 48.374 & 30.736 & 20.734 \\
		$2^{10}$ & 654.121 & 165.791 & 91.781 & 56.554 & 35.401 \\
		\bottomrule
	\end{tabular}
\end{table}

Meanwhile, by using the results of FSE-Shooting at $\omega_1 = 44.535$, $45.133$, $44.819$, $45.135$ as the initial conditions of VCHB constructed by the previous method $m$-VCF, the quasi-periodic solutions are computed as purple crosses in Fig. \ref{fig: FRCs d=e=2}a, where the parameters of $m$-VCF are defined as $K_2 = [0, 1, \dots, 4, 5]^\top$, $\mathrm{U}_2 = 11$ and $\mathrm{S}_1 = \mathrm{S}_2 = 2^5$. $K_1$ has three selections, namely $[0, 1, \dots, 5]^\top$ for points at $\omega_1 = 44.535$, $45.135$, $[0, 1, \dots, 7]^\top$ for $\omega_1 = 44.819$ and $[0, 1, \dots, 9]^\top$ for $\omega_1 = 45.133$. Correspondingly, there are $n\mathrm{U} = 28677$, $39105$ and $49533$ unknown coefficients $\mathbf{Q}^d$ in NAEs \eqref{EqA2}. And the computation costs per iteration of continuation of these three cases are respectively $61.141$s, $146.855$s and $283.112$s, which are significantly larger than that of the FSE-Shooting method with parallelization operations. In this example, the time interval of FSE-Shooting is between $7.709$s and $35.401$s.

\begin{figure}[ht] 
	\centering 
	\includegraphics[width=0.7\textwidth]{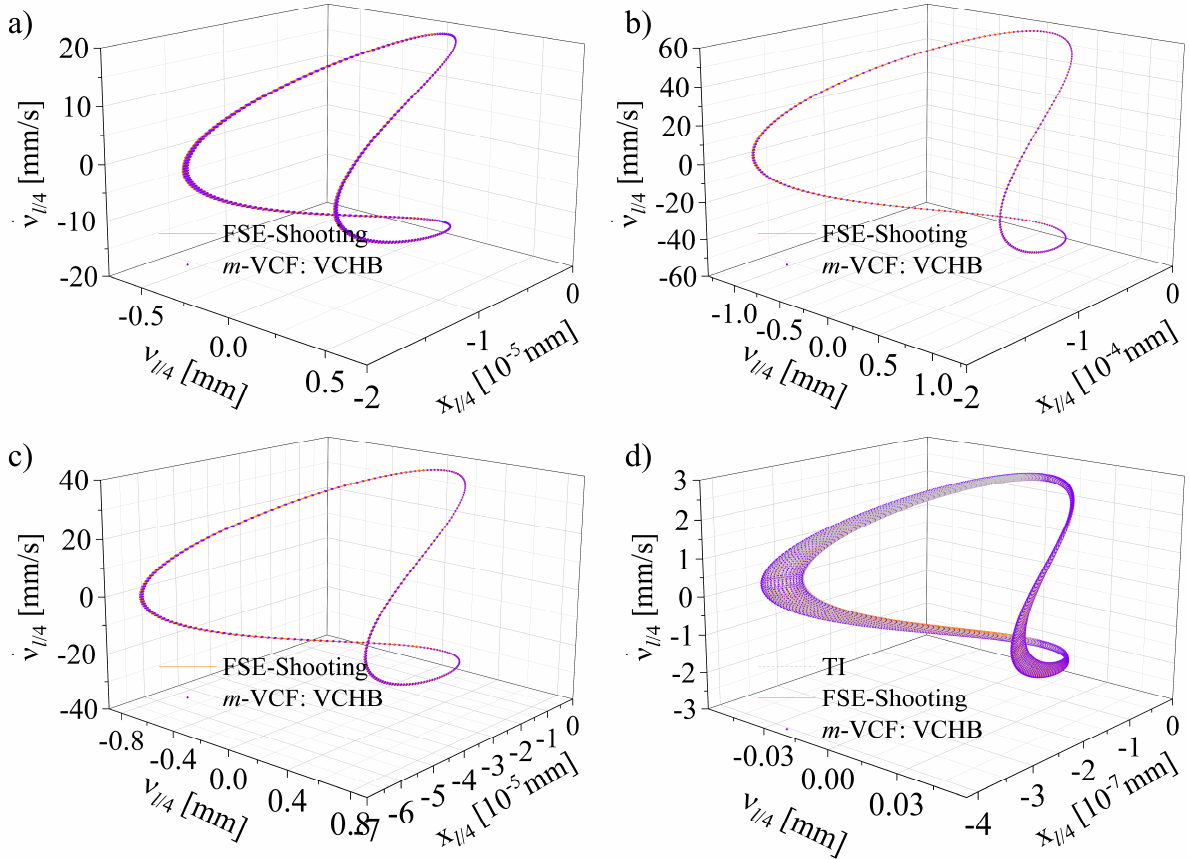} 
	\caption{\small The Phase diagrams of the clamped-clamped von Kármán beam under the quasi-periodic excitation with two base frequencies at: a): $\omega_1 = 44.535$, b): $\omega_1 = 45.133$, c): $\omega_1 = 44.819$, d): $\omega_1 = 45.135$.} 
	\label{fig: PhaseD d=e=2}
\end{figure}

\begin{figure}[htb] 
	\centering 
	\includegraphics[width=0.7\textwidth]{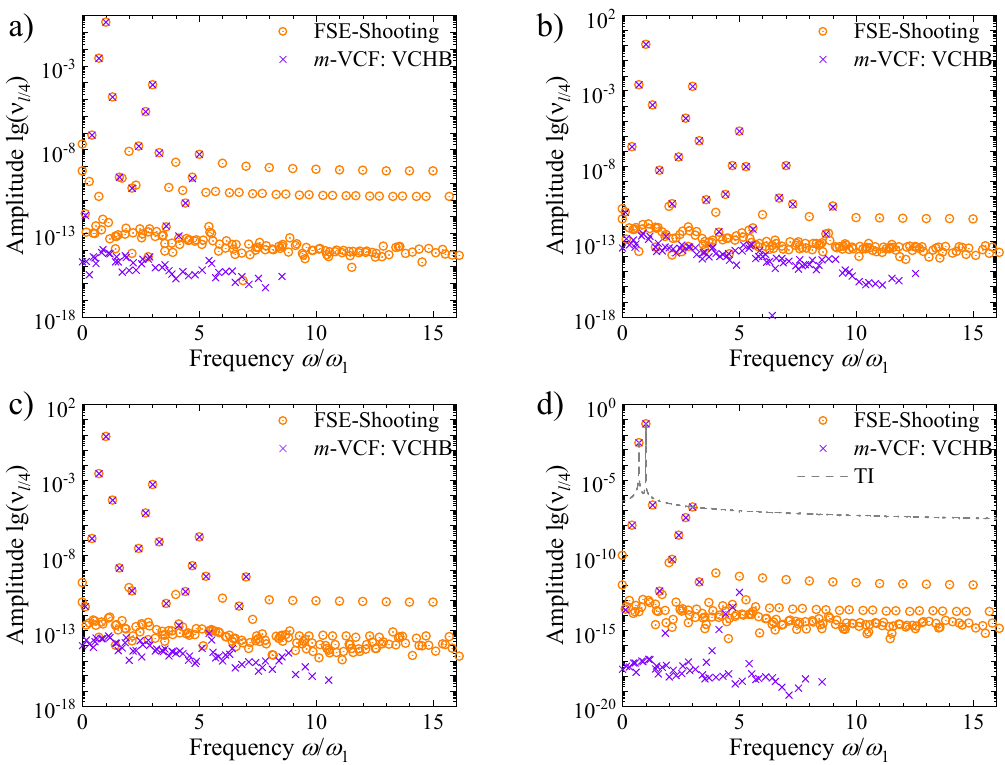} 
	\caption{\small The amplitude-frequency diagrams of the clamped-clamped von Kármán beam under the quasi-periodic excitation with two base frequencies at: a): $\omega_1 = 44.535$, b): $\omega_1 = 45.133$, c): $\omega_1 = 44.819$, d): $\omega_1 = 45.135$.} 
	\label{fig: AFD d=e=2}
\end{figure}

Fig. \ref{fig: PhaseD d=e=2} and Fig. \ref{fig: AFD d=e=2} show the Phase diagrams in the time domain and amplitude-frequency diagrams in the frequency domain at $\omega_1 = 44.535$, $45.133$, $44.819$, $45.135$. The results show that the solutions computed by these two methods have good consistency. In additions, the classical time integration method, (ode15s function of MATLAB) are used to compute quasi-periodic solutions at $\omega_1 = 45.135$. It can be found that, although its results are consistent with those of other methods in the time domain, it fails to reflect the complex characteristics of the quasi-periodic solutions in the frequency domain. And the time integration method is extremely time-consuming, taking about 21 hours.


\subsubsection{Quasi-periodic forcing with 3 base frequencies}\label{sec4.1.2}

Then, the case under the quasi-periodic excitation with $e = 3$ and quasi-periodic solutions with $d = 3$ is focused, the number of elements of beam is selected as $N_n = 12$. There are 33 DOFs in system \eqref{Eq25}, also 66 state variables in $1^\text{st}$ ODEs.

\begin{figure}[ht] 
	\centering 
	\includegraphics[width=0.45\textwidth]{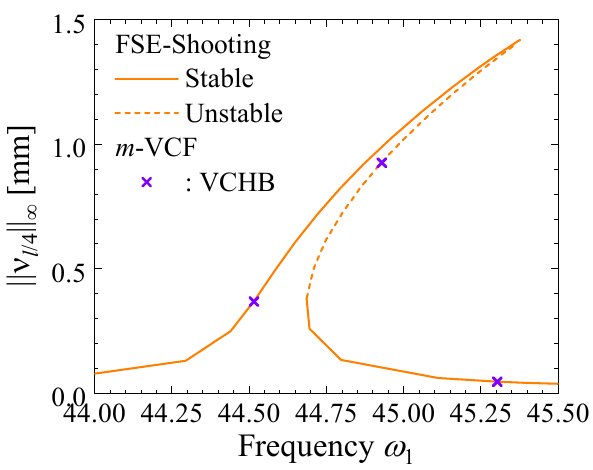} 
	\caption{\small The FRCs in the amplitude of transverse displacement at $l/4$ for clamped-clamped von Kármán beam under the quasi-periodic excitation with three base frequencies by using the proposed method (FSE-Shooting) and previous method (\textit{m}-VCF).} 
	\label{fig: FRCs d=e=3}
\end{figure}

The quasi-periodic solutions are computed by using the proposed method (FSE-Shooting) and previous method ($m$-VCF), as shown in Fig. \ref{fig: FRCs d=e=3}, where the solid curves and dotted curves are the stable quasi-periodic solutions and unstable ones computed by FSE-Shooting. Purple crosses are the results of VCHB at $\omega_1 = 44.516$, $44.930$, $45.303$, which is constructed by the $m$-VCF. The parameters of FSE-Shooting are defined as $K_2 = K_3 = [0, 1, \dots, 5]^\top$, $\mathrm{U}_2 = \mathrm{U}_3 = 11$, $\mathrm{S}_1 = 2^8$ and $\mathrm{S}_2 = \mathrm{S}_3 = 2^4$. And 32 Processors of CPU are used to compute the quasi-periodic solutions with three base frequencies. The parameters of the VCHB are $K_2 = K_3 = [0, 1, \dots, 5]^\top$, $\mathrm{U}_2 = \mathrm{U}_3 = 11$ and $\mathrm{S}_1 = \mathrm{S}_2 = \mathrm{S}_3 = 2^4$. Here, $K_1$ has two selections, namely $[0, 1, \dots, 5]^\top$ for points at $\omega_1 = 44.516$, $45.303$ and $[0, 1, \dots, 7]^\top$ for $\omega_1 = 44.930$. Therefore, there are 7986 unknown coefficients in FSE-Shooting and 43923 or 59895 unknown coefficients in VCHB. And the computation costs per iteration of continuation are $56.234$s in FSE-Shooting and $226.125$s or $469.667$s in VCHB. The proposed Shooting method is 4 times or even 8.4 times that of the previous discretization method.

\begin{figure}[htb] 
	\centering 
	\includegraphics[width=0.9\textwidth]{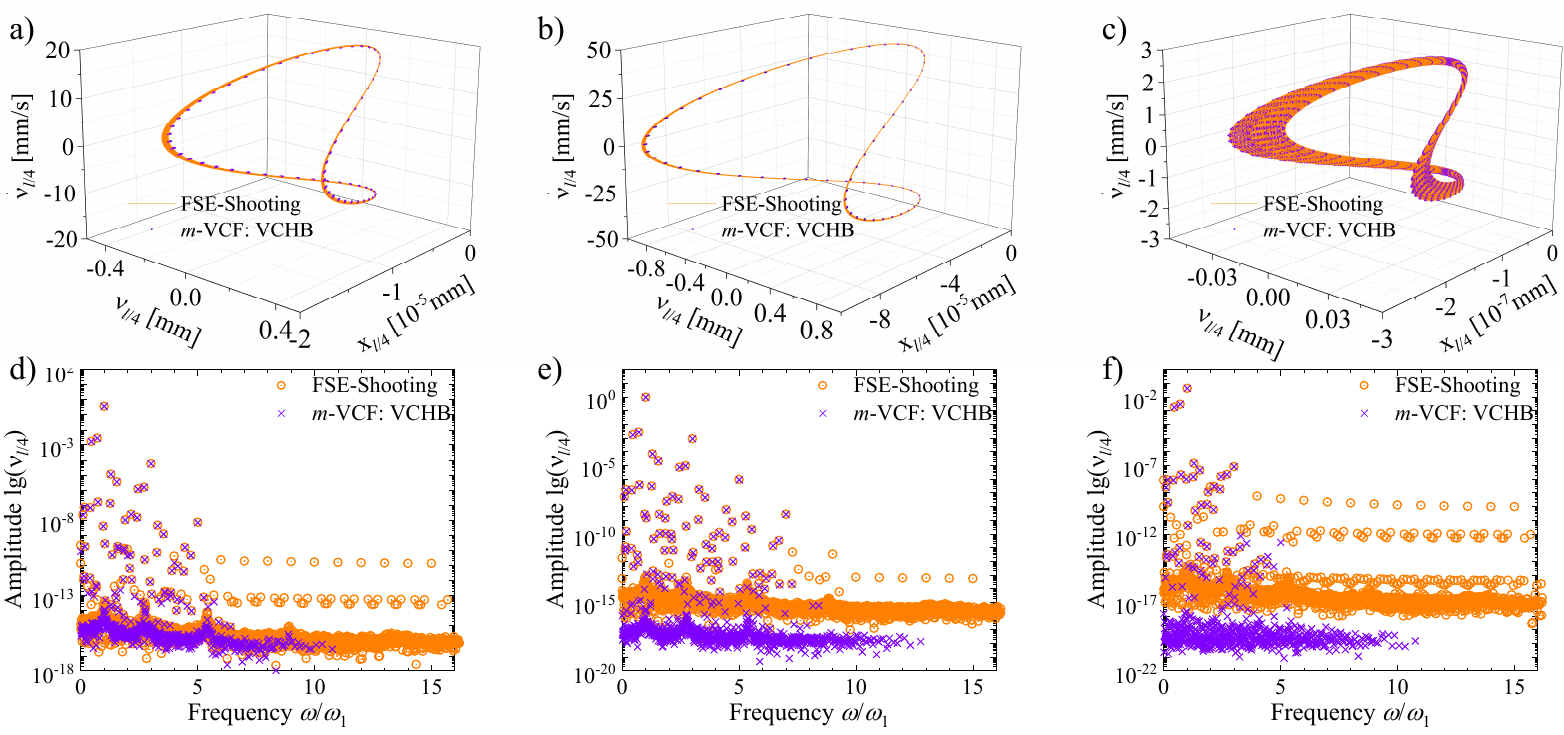} 
	\caption{\small The results of the clamped-clamped von Kármán beam under the quasi-periodic excitation with three base frequencies at $\omega_1 = 44.516$, $44.930$, $45.303$; a), b), and c): Phase diagrams; d), e) and f): Amplitude-frequency diagrams.} 
	\label{fig: PD & AFD d=e=3}
\end{figure}

\begin{figure}[htbp] 
	\centering 
	\includegraphics[width=0.7\textwidth]{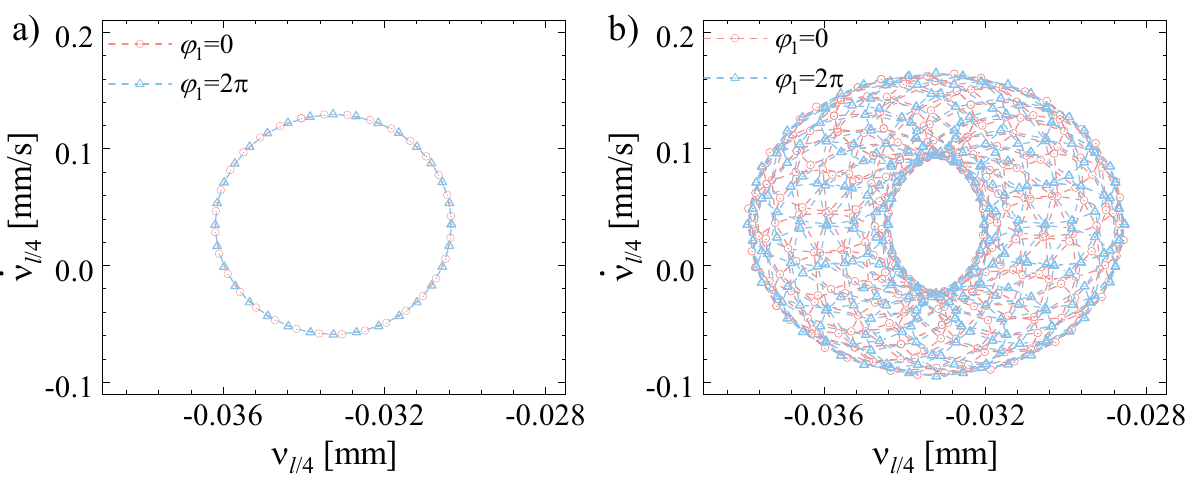} 
	\caption{\small The Poincaré sections of quasi-periodic solutions at $\varphi_1 = 0$, $2\pi$ under the quasi-periodic excitation at $\omega_1 = 45.5$, a): The limit cycles when $d = 2$; b): The limit torus when $d = 3$.} 
	\label{fig: initial & terminal states}
\end{figure}

Fig. \ref{fig: PD & AFD d=e=3} shows the Phase diagrams in the time domain and amplitude-frequency diagrams in the frequency domain of the clamped-clamped von Kármán beam under the quasi-periodic excitation with three base frequencies at $\omega_1 = 44.516$, $44.930$, $45.303$. The results also show the good consistency of the solutions computed by two methods. Compare these two cases with $d = 2, 3$, Fig. \eqref{fig: initial & terminal states} shows the initial (red) and terminal points (blue) of trajectories initialized at Poincaré sections under the quasi-periodic excitation at $\omega_1 = 45.5$, where the sections are respectively limit cycles for $d = 2$ and limit torus for $d = 3$. The aim of FSE-Shooting method is to find the initial conditions (Fourier coefficients) of $(d-1)$-tori at $\varphi_1 = 0$ by using the coupling conditions \eqref{EqCoupC}. In addition, as the number of base frequencies of excitation increases, the area of this Poincaré section increases accordingly, and the solutions of the system becomes more complex. 

\subsection{A periodically-forced von Kárman beam with Neimark-Sacker (NS) bifurcation}\label{sec4.2}

\begin{figure}[htb] 
	\centering 
	\includegraphics[width=0.75\textwidth]{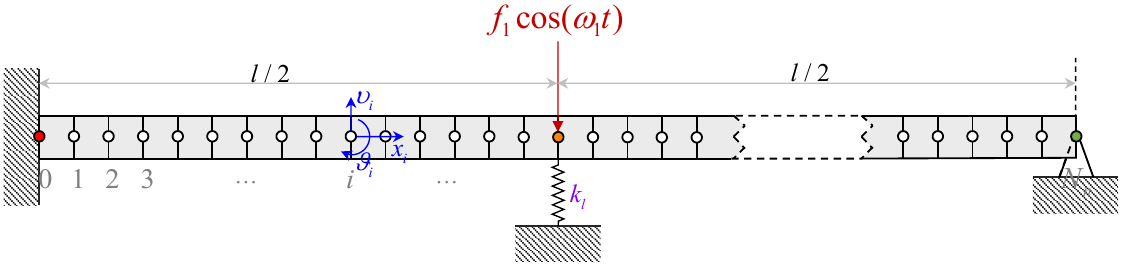} 
	\caption{\small A Finite-element clamped-pinned von Kármán beam with a support linear spring and periodic excitation at its midspan.} 
	\label{fig: CP-VKB}
\end{figure}

In Fig. \ref{fig: CP-VKB}, the second system is a clamped-pinned von Kármán beam with a support linear spring and periodic excitation at its midspan \cite{R70,R78}, whose quasi-periodic solutions occur after the Neimark-Sacker (NS) bifurcations of the periodic orbits \cite{R78}. Each node has three DOFs, namely, the axial displacement $x_i$, the transverse displacement $v_i$ and the rotation angle $\vartheta_i$. Hence, there are $3N_n - 2$ degrees of freedom, because $x_0, v_0, \vartheta_0, x_{N_n}, v_{N_n} = 0$. The Physical parameters of the von Kármán beam system are summarized in Table \ref{tab3}. Given the axial stretching due to large transverse deflections, the nonlinearity in this system is distributed. Due to the linear support spring at midspan, the system exhibits 1:3 internal resonance with natural frequency $\omega_{l_2} \approx 3\omega_{l_1}$, where the first two natural frequencies are $\omega_{l_1} = 33.20$ rad/s and $\omega_{l_2} = 99.59$ rad/s. Again, the system is modeled by finite-element model available in SSMTool.

\begin{table}[htbp]
	\centering
	\caption{Physical parameters of the clamped-pinned von Kármán beam system.}
	\begin{tabular}{ll}
		\toprule
		Physical parameter & Value \\
		\midrule
		Young's modulus $E$ [kPa] & $4.5 \times 10^7$ \\
		Density $\rho$ [kg/mm$^3$] & $1780 \times 10^{-9}$ \\
		Height of beam $h$ [mm] & 10 \\
		Width of beam $b$ [mm] & 10 \\
		Length of beam $l$ [mm] & 2700 \\
		Poisson's ratio $\nu$ & 0.3 \\
		Mass-proportional damping coefficient $\alpha$ [s] & 0 \\
		Stiffness-proportional damping coefficient $\beta$ [s$^{-1}$] & $2/9 \times 10^{-4}$ \\
		Stiffness of the linear support spring $k_l$ [N/mm] & 37 \\
		Amplitude of excitation $f_1$ [N] & 20 \\
		\bottomrule
	\end{tabular}
	\label{tab3}
\end{table}

\begin{figure}[htb] 
	\centering 
	\includegraphics[width=0.45\textwidth]{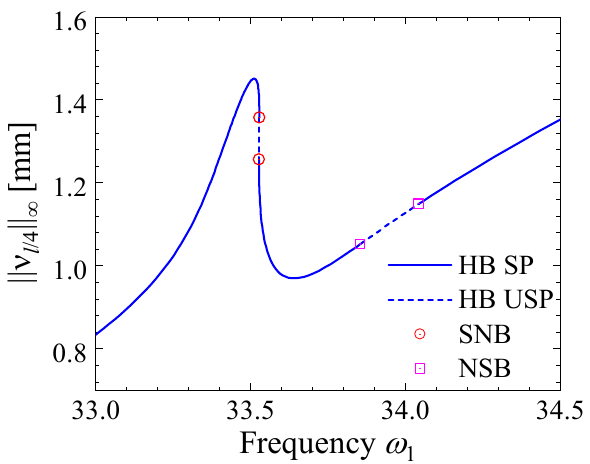} 
	\caption{\small The FRCs of periodic solutions in the amplitude of transverse displacement at $l/4$ for clamped-pinned von Kármán beam with 94-DOFs by using the HB method.} 
	\label{fig: FRCs d=e=1}
\end{figure}

The Fig. \ref{fig: FRCs d=e=1} shows the FRCs of periodic solutions of the von Kármán beam with 94 DOFs (32 elements), where the solid curves represent stable solutions and the dotted curves stand for unstable solutions, the squares represent the Neimark-Sacker bifurcations (NSBs) and circles are saddle-node bifurcations (SNBs). After the NSBs, the unstable periodic solutions are easily perturbed into the stable quasi-periodic solutions, where the second base frequency of solution are a priori unknown, i.e., $d = 2$, $e = 1$.

\begin{figure}[htb] 
	\centering 
	\includegraphics[width=0.8\textwidth]{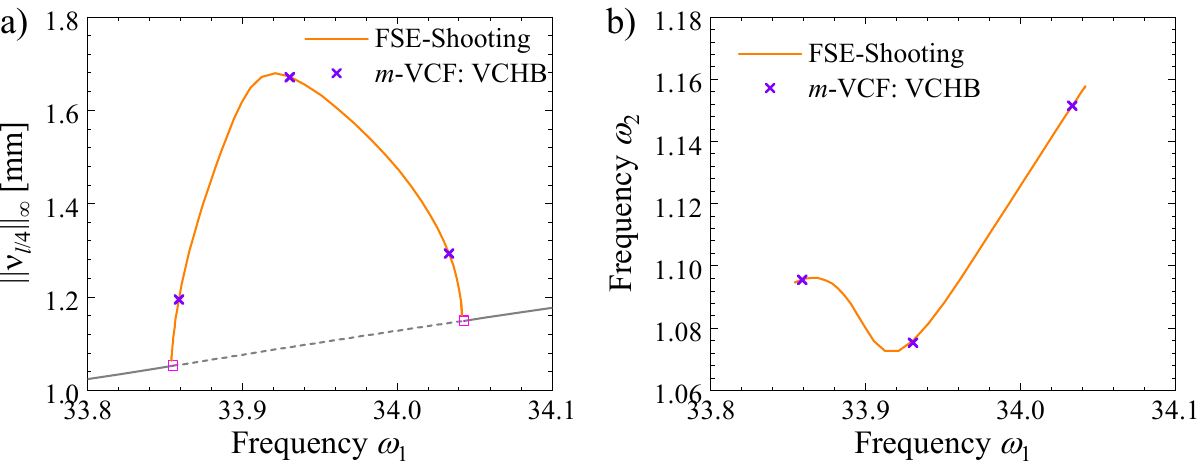} 
	\caption{\small The results for clamped-pinned von Kármán beam with 94-DOFs by using the proposed method (FSE-Shooting) and previous method ($m$-VCF), a): FRCs of quasi-periodic solutions in the amplitude of transverse displacement at $l/4$; b): Relationship between the base frequencies of the quasi-periodic solutions.
	} 
	\label{fig: FRCs d=2 e=1}
\end{figure}

Based on the initialization technique in \cite{R41,R50} and the phase conditions introduced in Sub-section \ref{sec3.2.2}, the quasi-periodic solutions tracked by the FSE-Shooting and previous method ($m$-VCF) are depicted in Fig. \ref{fig: FRCs d=2 e=1}a, where the solutions are all stable. Here, the parameters of the FSE-Shooting are $K_2 = [0, 1, \dots, 7]^\top$, $\mathrm{U}_2 = 15$, $\mathrm{S}_1 = 2^{10}$ and $\mathrm{S}_2 = 2^5$. Therefore, there are 2820 unknown coefficients $\mathbf{Z}(0)$ in NAEs \eqref{EqNAE_FSES}. By using 32 Processors of CPU, the time per iteration is about $12.322$s in the continuation. By using the results of FSE-Shooting at $\omega_1 = 33.859$, $33.931$, $34.034$ as the initial conditions of VCHB constructed by the previous method $m$-VCF, the quasi-periodic solutions are computed as purple crosses in Fig. \ref{fig: FRCs d=2 e=1}a, whose parameters are defined as $K_2 = [0, 1, \dots, 7]^\top$, $\mathrm{U}_2 = 15$ and $\mathrm{S}_1 = \mathrm{S}_2 = 2^5$. $K_1$ has two selections, namely $[0, 1, \dots, 11]^\top$ for points at $\omega_1 = 33.859$, $34.034$ and $[0, 1, \dots, 13]^\top$ for $\omega_1 = 33.931$. Therefore, there are 32430 or 38070 unknown coefficients. And the computation costs per iteration of continuation are $69.985$s and $113.885$s, respectively. The proposed shooting method is 5.8 times or even 9.5 times that of the previous discretization method.

\begin{figure}[htb] 
	\centering 
	\includegraphics[width=0.9\textwidth]{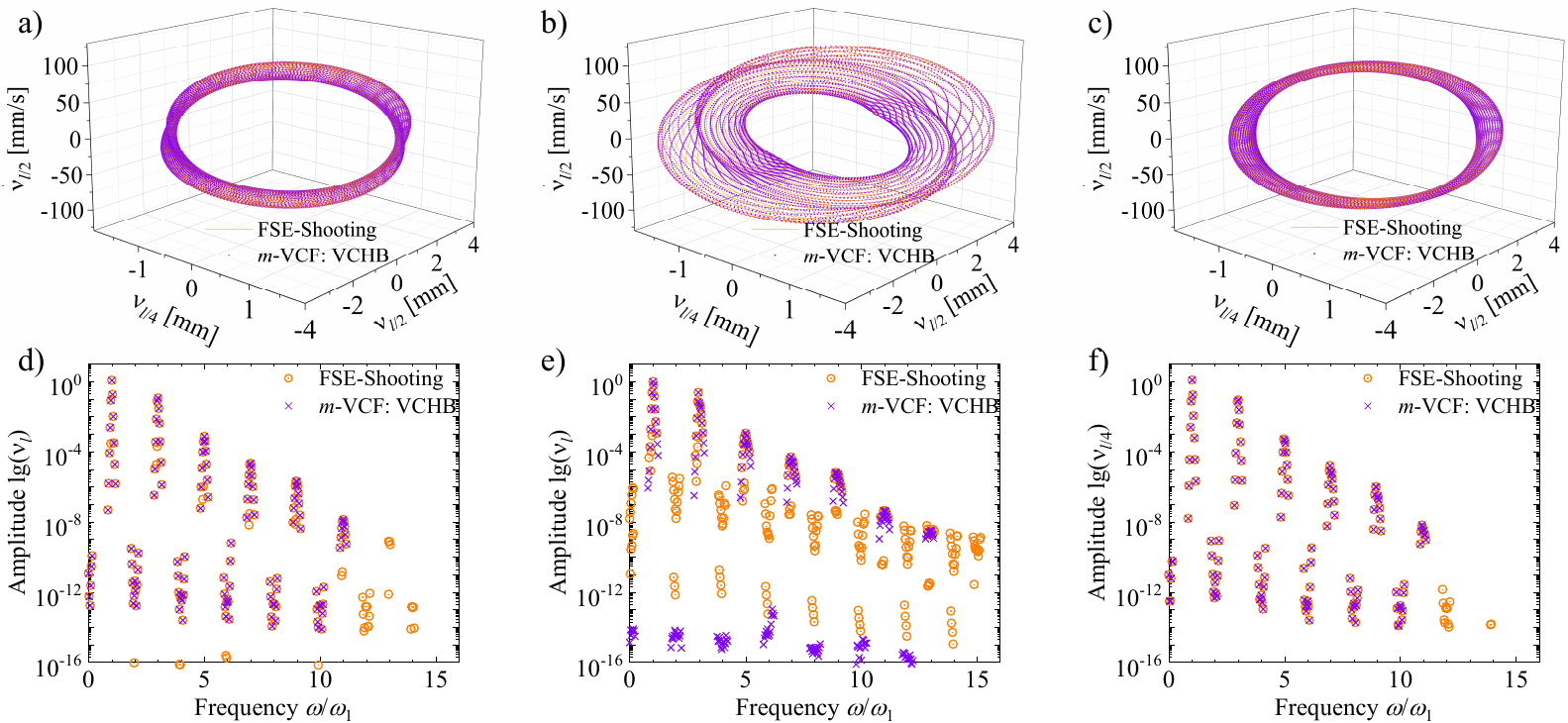} 
	\caption{\small The results of quasi-periodic solutions of the clamped-pinned von Kármán beam under the periodic excitation at $\omega_1 = 33.859$, $33.931$, $34.034$; a), b), and c): Phase diagrams; d), e) and f): Amplitude-frequency diagrams.} 
	\label{fig: PD & AFD d=2 e=1}
\end{figure}

More importantly, the relationship of base frequencies of quasi-periodic solutions are shown in Fig. \ref{fig: FRCs d=2 e=1}b, where the curves are the results of FSE-Shooting method and crosses are those of previous method ($m$-VCF). Fig. \ref{fig: PD & AFD d=2 e=1} shows the Phase diagrams in the time domain and amplitude-frequency diagrams in the frequency domain of quasi-periodic solutions of the clamped-pinned von Kármán beam under the periodic excitation. These results all show the good consistency of the solutions computed by two methods. Among them, the proposed FSE-Shooting method can provide much richer dynamical behaviors in the frequency domain.

\subsection{A quasi-periodically-forced shallow-arch shell structure}\label{sec4.3}

\begin{figure}[htb] 
	\centering 
	\includegraphics[width=0.85\textwidth]{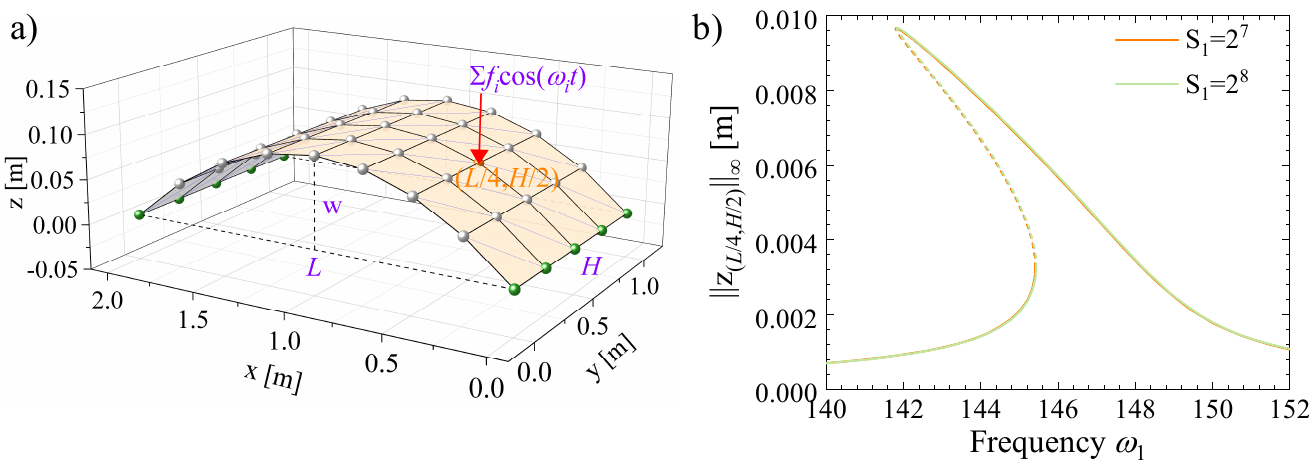} 
	\caption{\small a): A Finite-element shallow-arch structure with $N_H \times N_L$ meshes; b): The FRCs of quasi-periodic solutions of shallow-arch shell structure by using the proposed method (FSE-Shooting).} 
	\label{fig: shallow-arch shell}
\end{figure}

In this Sub-section, a quasi-periodically-forced shallow-arch shell structure is adapted \cite{R18,R70}, as illustrated in Fig. \ref{fig: shallow-arch shell}a, where the shell is simply supported at the two opposite edges aligned along the y-axis. The Physical parameters of this shell are given in Table \ref{tab4}. The curvature parameter $w$ is defined as the height of the midpoint relative to the end. The quasi-periodic excitation $\sum_{i=1}^d f_i \cos(\omega_i t)$ is forced at the position $(L/4, H/2)$ along the z-axis, and the base frequencies have the rule of $\rho_2 = \omega_2 / \omega_1 = 1/\sqrt{2}$. Using flat, triangular shell elements, $N_H \times N_L$ ($N_L = 2N_H$) meshes are used to discrete this model. The discrete model here contains $2N_H N_L$ elements, and each node in the elements has six DOFs, namely the displacements $\mathrm{x}, \mathrm{y}, \mathrm{z}$ and the angle displacements $\vartheta_{\mathrm{x}}, \vartheta_{\mathrm{y}}, \vartheta_{\mathrm{z}}$. $d = 2$ and $e = 2$ hold in this example.

\begin{table}[htbp]
	\centering
	\caption{Physical parameters of the shallow-arch shell structure system.}
	\begin{tabular}{ll}
		\toprule
		Physical parameter & Value \\
		\midrule
		Young's modulus $E$ [GPa] & 70 \\
		Density $\rho$ [kg/m$^3$] & 2700 \\
		Length of shell $L$ [m] & 2 \\
		Width of shell $H$ [m] & 1 \\
		Thickness of shell $t$ [m] & 0.01 \\
		Curvature parameter $w$ [m] & 0.1 \\
		Poisson's ratio $\nu$ & 0.33 \\
		Viscous damping rate of material [Pa s] & $10^5$ \\
		Amplitude of excitation $f_1$ and $f_2$ [N] & $40, 20$ \\
		\bottomrule
	\end{tabular}
	\label{tab4}
\end{table}

First, the shallow-arch structure is discretized into 64 elements by using $4 \times 8$ meshes, where the 240 DOFs are introduced in the finite-element shell system. The FRCs of quasi-periodic solutions $\mathrm{z}_{(L/4, H/2)}$ computed by FSE-Shooting method with the step parameter $\mathrm{S}_1 = 2^7$, $2^8$ are compared in Fig. \ref{fig: shallow-arch shell}b, where the slide curves denote the stable quasi-periodic solutions and the dotted curves are unstable ones. Other parameters are defined as $K_2 = [0, 1, \dots, 5]^\top$, $\mathrm{U}_2 = 11$ and $\mathrm{S}_2 = 2^5$. The number of unknown coefficients $\mathbf{Z}(0)$ in NAEs \eqref{EqNAE_FSES} is 5280. Results show that the FRCs exhibit good consistency for two $\mathrm{S}_1 = 2^7$, $2^8$. And the computation costs per iteration of continuation are respectively $21.658$s and $45.505$s. Considering the computational cost, $\mathrm{S}_1 = 2^7$ is used in the following numerical example.

\begin{figure}[htb] 
	\centering 
	\includegraphics[width=0.8\textwidth]{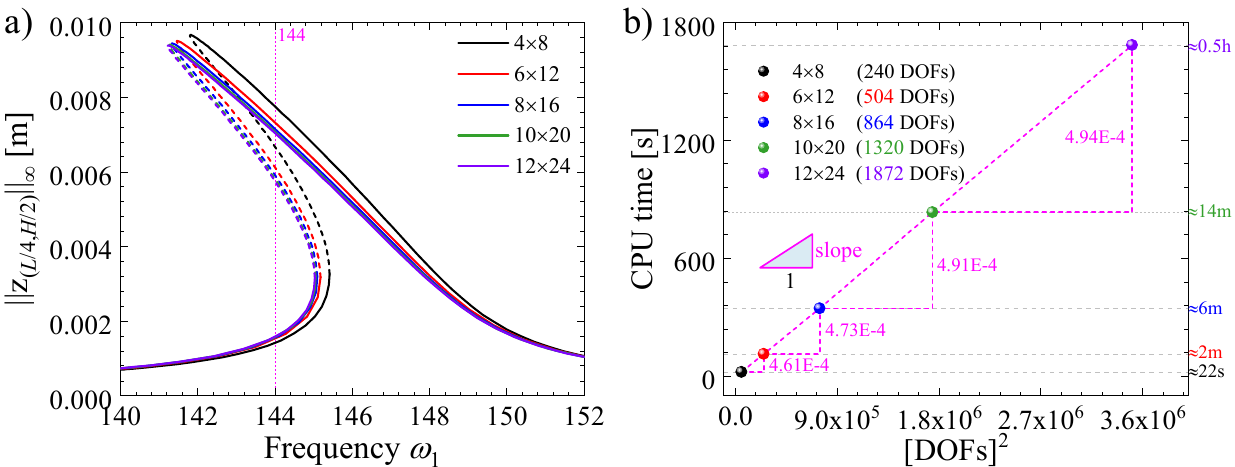} 
	\caption{\small a): The FRCs of quasi-periodic solutions $\mathrm{z}_{(L/4, H/2)}$ computed by the FSE-Shooting method by using different meshes; b) The computation cost per iteration in the continuation.} 
	\label{fig: shallow-arch shell FRCs d=e=2}
\end{figure}

Then, by introducing $4 \times 8$, $6 \times 12$, $8 \times 16$, $10 \times 20$ and $12 \times 24$ meshes, the shell is modeled as the Finite-element systems with 240, 504, 864, 1320 and 1872 DOFs. Here, the numbers of the Fourier coefficients in NAEs \eqref{EqNAE_FSES} are respectively 5280, 11088, 19008, 29040 and 41184. The FRCs of quasi-periodic solution $\mathrm{z}_{(L/4, H/2)}$ are computed in Fig. \ref{fig: shallow-arch shell FRCs d=e=2}a, where the slide curves represent stable solutions and dotted curves are unstable ones. As the mesh is refined, the FRCs gradually stabilize, and the numerical solutions converge to the true solutions. The relationship between the computation cost per iteration in the continuation by using 32 workers of CPU and the square of degrees of freedom is illustrated in Fig. \ref{fig: shallow-arch shell FRCs d=e=2}b, and they are nearly in a proportional relationship with similar slope.

\begin{figure}[htb] 
	\centering 
	\includegraphics[width=0.85\textwidth]{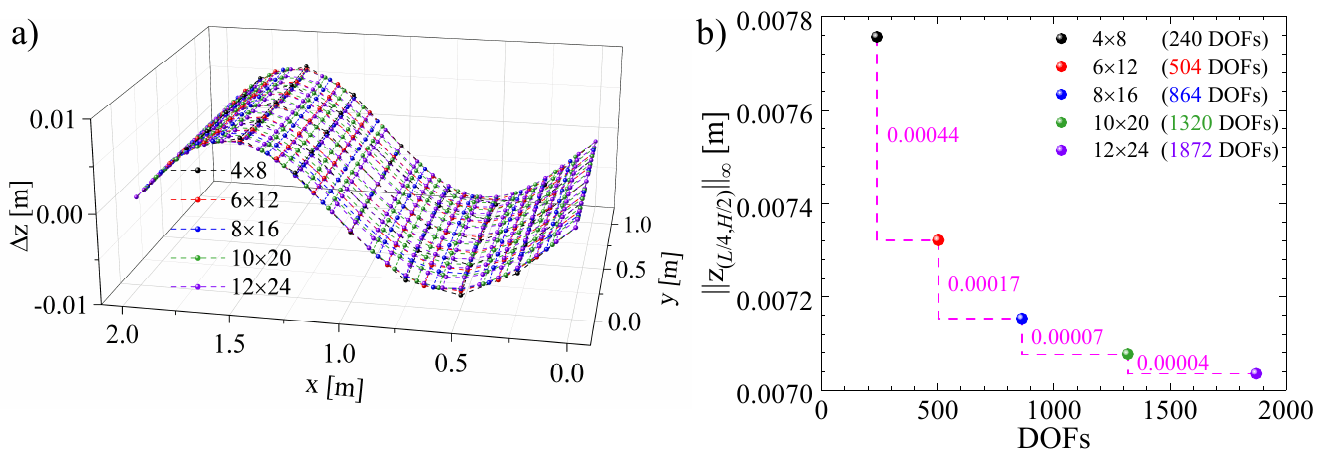} 
	\caption{\small a): The deformation of the shell along the z-axis when the absolute value of $\mathrm{z}_{(L/4, H/2)}$ is largest; b) The differences between the largest absolute value of $\mathrm{z}_{(L/4, H/2)}$ with different meshes.} 
	\label{fig: shallow-arch shell -Time}
\end{figure}

\begin{figure}[htb] 
	\centering 
	\includegraphics[width=0.95\textwidth]{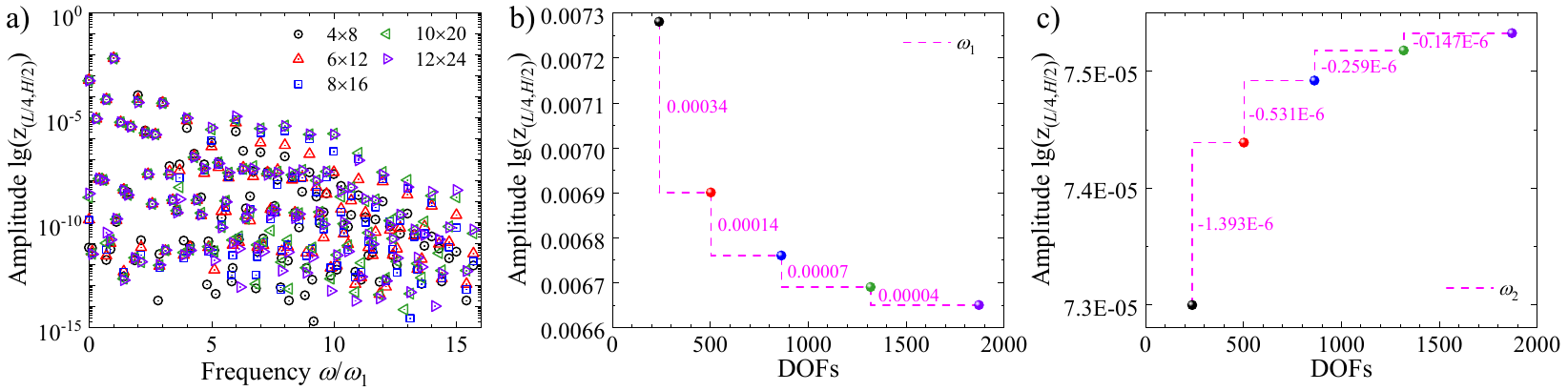} 
	\caption{\small a): The differences of amplitude of $\mathrm{z}_{(L/4, H/2)}$ in the frequency domain; b): The amplitude of frequency component $\omega = \omega_1$; c): The amplitude of frequency component $\omega = \omega_2$.} 
	\label{fig: shallow-arch shell -Frequency}
\end{figure}

Further, Fig. \ref{fig: shallow-arch shell -Time}a shows the convergence of the deformation of shell at $\omega_1 = 144$ along the z-axis when the absolute value of $\mathrm{z}_{(L/4, H/2)}$ is largest. The differences between the largest absolute value of $\mathrm{z}_{(L/4, H/2)}$ as the mesh is refined are illustrated in Fig. \ref{fig: shallow-arch shell -Time}b, which is increasingly small. The same results can be concluded in Fig. \ref{fig: shallow-arch shell -Frequency}, which shows the differences of amplitude of $\mathrm{z}_{(L/4, H/2)}$ in the frequency domain. This indicates that the finite-element model has good convergence, and the computational results are reliable when the mesh is sufficiently dense.

\section{Conclusions}\label{sec5}

In this work, a Fourier series expansion-based shooting method in the multi-frequency domain (FSE-Shooting) is proposed for the parallel computation of quasi-periodic solutions containing $d$ base frequencies in finite element problems, where the base frequencies can be unknown. By defining a new hyper-time domain via multi-period parameterization, the $d$-torus is parameterized as a collection of trajectories initialized at its first-order Poincaré section under this new multi-trajectory parameterization. By means of Fourier series expansion, we analytically derive the nonlinear algebraic equations (NAEs) for searching the Fourier coefficients of the initial state of the trajectories, where AFT and time integration are used to compute Fourier coefficients of the terminal state in parallel. In addition, a robust phase condition is also proposed to fix the phase of $d$-torus with respect to unknown base frequencies. Numerical examples demonstrate its high efficiency potential in high-dimensional systems, including a shell model with 1872 DOFs in Sub-section \ref{sec4.3}. In addition, like the classical shooting method for periodic solutions, the by-product of FSE-Shooting can be also directly used to assess the stabilities of quasi-periodic solutions in terms of Lyapunov exponents, as shown in \ref{appF}.

 Specially, the results of Example in Sub-section \ref{sec4.1.1} show that the efficiency is greater than 0.87 when using 4 Processors. Although it decreases to some extent as the number of Processors increases, the speedup of this parallelized computation nearly reaches a maximum of 18.5 times when using 32 Processors. Actually, the reduction of efficiency should be attributed to the increased difficulty in calculator resource and worker redundancy. In addition, this work only utilizes the "parfor" function of MATLAB for simple CPU parallelization, and the parallel efficiency of this method still has room for improvement. By the way, the potential of its parallelized computation may become even more prominent as $d$ increases. For example, if the computer is equipped with 256 processors, the time cost of example in Sub-section \ref{sec4.1.2} will be reduced more significantly than when using 32 processors. In a word, the results show that the proposed shooting method is about $7$-$12$ times more efficient than the best-performing discretization method.
 
Moreover, we provide a simple self-excited system in \ref{appG} to illustrate the applicability of FSE-Shooting method when multiple base frequencies are unknown, where the quasi-periodic solution of the system is generated by secondary Hopf bifurcation. In particular, we find that discretization methods may be a better choice for low-dimensional systems due to their simpler structure. Also, there are some future recommendations for the FSE-Shooting: 1) a multiple FSE-Shooting method is recommended to improve the robustness and convergence domain; 2) the analytical expression of Jacobian matrix with respect to $p$ can be derived to replace the finite difference procedure. 3) the extension of FSE-Shooting to nonsmooth dynamical systems.

\section*{Acknowledgments}

This work was supported by the National Natural Science Foundation of China (Grant No. 12172267 and 12302014).

\appendix
\section{The discretization techniques}\label{appA}

There are two types of discretization process: namely, full-discretization (single-step discretization) and semi-discretization (multi-step discretization) methods. Before introducing these two methods, the Eq. \eqref{2ODEs} is rewritten as:
\begin{equation}
	\begin{aligned}
	\mathbf{R}^{0} (\mathbf{Q}^{0} (\boldsymbol{\tau}_1^d), \boldsymbol{\tau}_1^d) & = \mathbf{M}^0\ddot{\mathbf{Q}}^0(\boldsymbol{\tau}_1^d)+\mathbf{D}^0\dot{\mathbf{Q}}^0(\boldsymbol{\tau}_1^d)+\mathbf{K}^0\mathbf{Q}^0(\boldsymbol{\tau}_1^d) \\& +\mathbf{\Theta}^0[\mathbf{F}^0(\mathbf{Q}^0,\dot{\mathbf{Q}}^0)-\mathbf{E}^0(\boldsymbol{\tau}_1^e)]=\mathbf{0},\label{EqA0}
	\end{aligned}
\end{equation}
For full-discretization methods, a single-step constant-coefficients formulation ($s$-CCF) is proposed to approximate the quasi-periodic solution and its derivatives:
\begin{equation}
	\begin{aligned}
		\mathbf{Q}^0(\boldsymbol{\tau}_1^{d}) &= \left( \mathbf{I}_n \otimes \boldsymbol{\Phi}(\boldsymbol{k}_1^{d}, \boldsymbol{\tau}_1^{d}) \right) \mathbf{Q}^d(\boldsymbol{k}_1^{d}); \\
		\dot{\mathbf{Q}}^0(\boldsymbol{\tau}_1^{d}) &= \left( \mathbf{I}_n \otimes \dot{\boldsymbol{\Phi}}(\boldsymbol{k}_1^{d}, \boldsymbol{\tau}_1^{d}) \right) \mathbf{Q}^d(\boldsymbol{k}_1^{d}) \\
		&= \left( \mathbf{I}_n \otimes \sum_{i=1}^{d} \omega_i \partial_{\boldsymbol{\tau}_i} \boldsymbol{\Phi}(\boldsymbol{k}_1^{d}, \boldsymbol{\tau}_1^{d}) \right) \mathbf{Q}^d(\boldsymbol{k}_1^{d}); \\
		\ddot{\mathbf{Q}}^0(\boldsymbol{\tau}_1^{d}) &= \left( \mathbf{I}_n \otimes \ddot{\boldsymbol{\Phi}}(\boldsymbol{k}_1^{d}, \boldsymbol{\tau}_1^{d}) \right) \mathbf{Q}^d(\boldsymbol{k}_1^{d}) \\
		&= \left( \mathbf{I}_n \otimes \sum_{j=1}^{d} \sum_{i=1}^{d} \omega_j \omega_i \partial_{\boldsymbol{\tau}_j \boldsymbol{\tau}_i}^2 \boldsymbol{\Phi}(\boldsymbol{k}_1^{d}, \boldsymbol{\tau}_1^{d}) \right) \mathbf{Q}^d(\boldsymbol{k}_1^{d}),
	\end{aligned}
	\label{EqA1}
\end{equation}
where $\otimes$ is the Kronecker tensor product. $\mathbf{I}_n$ is $n \times n$ identity matrix. The quasi-periodic shape functions $\boldsymbol{\Phi}(\boldsymbol{k}_1^{d}, \boldsymbol{\tau}_1^{d}): \mathbb{U}^d \times \mathbb{T}^d \mapsto \mathbb{R}^{1 \times \mathrm{U}}$ and their derivatives are constructed by the $\boldsymbol{k}_1^{d}$-parameters $\boldsymbol{k}_1^{d} \in \mathbb{U}^d$, with $\mathbb{U}^d$ being a $d$-dimensional discrete domain. $\partial_{\tau_i} \boldsymbol{\Phi}(\boldsymbol{k}_1^{d}, \boldsymbol{\tau}_1^{d})$ and $\partial_{\tau_j \tau_i}^2 \boldsymbol{\Phi}(\boldsymbol{k}_1^{d}, \boldsymbol{\tau}_1^{d})$ are the derivatives of shape functions with respect to $\tau_i$ and $\tau_j \tau_i$, respectively. Substituting \eqref{EqA1} into Eq. \eqref{2ODEs} and discretizing $\boldsymbol{\Phi}(\boldsymbol{k}_1^{d}, \boldsymbol{\tau}_1^{d})$ and their derivatives, yields a set of nonlinear algebraic equations (NAEs) with respect to constant coefficients $\mathbf{Q}^d(\boldsymbol{k}_1^{d})$ and base frequencies $\boldsymbol{\omega}_1^d$:
\begin{equation}
	\mathbf{R}^d (\mathbf{Q}^d (\boldsymbol{k}_1^d), \boldsymbol{\omega}_1^d) = \mathbf{K}^d (\boldsymbol{\omega}_1^d) \mathbf{Q}^d (\boldsymbol{k}_1^d) + \boldsymbol{\Theta}^d (\boldsymbol{k}_1^d) \times [\mathbf{F}^d (\mathbf{Q}^d) - \mathbf{E}^d] = \mathbf{0},
	\label{EqA2}
\end{equation}
with
\begin{equation}
	\begin{aligned}
		{\mathbf{K}}^d &= \mathbf{M}^0 \otimes \tilde{\boldsymbol{\Upsilon}}^2 + \mathbf{D}^0 \otimes \tilde{\boldsymbol{\Upsilon}}^1 + \mathbf{K}^0 \otimes \tilde{\boldsymbol{\Upsilon}}^0; \quad {\boldsymbol{\Theta}}^d = \boldsymbol{\Theta}^0 \otimes \tilde{\boldsymbol{\Upsilon}}^0 \\
		\tilde{\boldsymbol{\Upsilon}}^0 &= \boldsymbol{\Upsilon}^0; \quad \tilde{\boldsymbol{\Upsilon}}^1 = \sum_{i=1}^{d} \omega_i \boldsymbol{\Upsilon}^{1,i}; \quad \tilde{\boldsymbol{\Upsilon}}^2 = \tilde{\boldsymbol{\Upsilon}}^1 \tilde{\boldsymbol{\Upsilon}}^1,
	\end{aligned}
	\label{EqA3}
\end{equation}
Here, $\boldsymbol{\Upsilon}^0$ and $\boldsymbol{\Upsilon}^{1,i}$ are $\mathrm{U} \times \mathrm{U}$ constant matrices, which are discretized form $\boldsymbol{\Phi}(\boldsymbol{k}_1^d, \boldsymbol{\tau}_1^d)$ and $\partial_{\tau_i} \boldsymbol{\Phi}(\boldsymbol{k}_1^d, \boldsymbol{\tau}_1^d)$.

For semi-discretization methods, a multi-step variable-coefficients formulation ($m$-VCF) is recently proposed in \cite{R41}. In $m$-VCF, the discretization system in $i$-1 step is governed by:
\begin{equation}
	\begin{aligned}
		\mathbf{R}^{i-1} (\mathbf{Q}^{i-1} (\boldsymbol{k}_1^{i-1}, \boldsymbol{\tau}_i^d), \boldsymbol{\omega}_1^{i-1}, \boldsymbol{\tau}_i^d) &= \mathbf{M}^{i-1} (\boldsymbol{\omega}_1^{i-1}) \ddot{\mathbf{Q}}^{i-1} (\boldsymbol{k}_1^{i-1}, \boldsymbol{\tau}_i^d) \\
		&+ \mathbf{D}^{i-1} (\boldsymbol{\omega}_1^{i-1}) \dot{\mathbf{Q}}^{i-1} (\boldsymbol{k}_1^{i-1}, \boldsymbol{\tau}_i^d) \\
		&+ \mathbf{K}^{i-1} (\boldsymbol{\omega}_1^{i-1}) \mathbf{Q}^{i-1} (\boldsymbol{k}_1^{i-1}, \boldsymbol{\tau}_i^d) \\
		&+ \boldsymbol{\Theta}^{i-1} \times [\mathbf{F}^{i-1} (\mathbf{Q}^{i-1}, \dot{\mathbf{Q}}^{i-1}) - \mathbf{E}^{i-1} (\boldsymbol{k}_1^{i-1}, \boldsymbol{\tau}_i^e)] \\
		&= \mathbf{0},
	\end{aligned}
	\label{EqA4}
\end{equation}
By approximating $\mathbf{Q}^{i-1}$ and its derivatives
\begin{equation}
	\begin{aligned}
		\mathbf{Q}^{i-1} (\boldsymbol{k}_1^{i-1}, \boldsymbol{\tau}_i^d) &\approx [\mathbf{I}_{n \mathrm{U}_1^{i-1}} \otimes \boldsymbol{\Phi}_i (k_i, \tau_i)] \mathbf{Q}^i (\boldsymbol{k}_1^i, \boldsymbol{\tau}_{i+1}^d); \\
		\dot{\mathbf{Q}}^{i-1} (\boldsymbol{k}_1^{i-1}, \boldsymbol{\tau}_i^d) &\approx [\mathbf{I}_{n \mathrm{U}_1^{i-1}} \otimes \omega_i {\boldsymbol{\Phi}}^{\prime}_i (k_i, \tau_i)] \mathbf{Q}^i (\boldsymbol{k}_1^i, \boldsymbol{\tau}_{i+1}^d) \\
		&\quad + [\mathbf{I}_{n \mathrm{U}_1^{i-1}} \otimes \boldsymbol{\Phi}_i (k_i, \tau_i)] \dot{\mathbf{Q}}^i (\boldsymbol{k}_1^i, \boldsymbol{\tau}_{i+1}^d); \\
		\ddot{\mathbf{Q}}^{i-1} (\boldsymbol{k}_1^{i-1}, \boldsymbol{\tau}_i^d) &\approx [\mathbf{I}_{n \mathrm{U}_1^{i-1}} \otimes \omega_i^2 {\boldsymbol{\Phi}}^{\prime\prime}_i (k_i, \tau_i)] \mathbf{Q}^i (\boldsymbol{k}_1^i, \boldsymbol{\tau}_{i+1}^d) \\
		&\quad + [\mathbf{I}_{n \mathrm{U}_1^{i-1}} \otimes 2 \omega_i {\boldsymbol{\Phi}}^{\prime}_i (k_i, \tau_i)] \dot{\mathbf{Q}}^i (\boldsymbol{k}_1^i, \boldsymbol{\tau}_{i+1}^d) \\
		&\quad + [\mathbf{I}_{n \mathrm{U}_1^{i-1}} \otimes \boldsymbol{\Phi}_i (k_i, \tau_i)] \ddot{\mathbf{Q}}^i (\boldsymbol{k}_1^i, \boldsymbol{\tau}_{i+1}^d),
	\end{aligned}
	\label{EqA5}
\end{equation}
where $\boldsymbol{\Phi}_i (k_i, \tau_i): \mathbb{U} \times \mathbb{T} \mapsto \mathbb{R}^{1 \times \mathrm{U}_i}$ represents a set of periodic functions with respect to $\boldsymbol{\tau}_i$. By submitting Eq. \eqref{EqA5} into Eq. \eqref{EqA4}, it is transformed into the differential system of $\mathbf{R}^i (\mathbf{Q}^i (\boldsymbol{k}_1^i, \boldsymbol{\tau}_{i+1}^d), \boldsymbol{\omega}_1^i, \boldsymbol{\tau}_{i+1}^d)$ with
\begin{equation}
	\begin{aligned}
		\mathbf{M}^i &= \mathbf{M}^{i-1} \otimes \boldsymbol{\Upsilon}_i^0; \\
		\mathbf{D}^i &= 2 \mathbf{M}^{i-1} \otimes \omega_i \boldsymbol{\Upsilon}_i^1 + \mathbf{D}^{i-1} \otimes \boldsymbol{\Upsilon}_i^0; \\
		\mathbf{K}^i &= \mathbf{M}^{i-1} \otimes \omega_i^2 \boldsymbol{\Upsilon}_i^2 + \mathbf{D}^{i-1} \otimes \omega_i \boldsymbol{\Upsilon}_i^1 + \mathbf{K}^{i-1} \otimes \boldsymbol{\Upsilon}_i^0; \\
		\boldsymbol{\Theta}^i &= \boldsymbol{\Theta}^{i-1} \otimes \boldsymbol{\Upsilon}_i^0,
	\end{aligned}
	\label{EqA6}
\end{equation}
where $\boldsymbol{\Upsilon}_i^0 ({k}_i)$, $\boldsymbol{\Upsilon}_i^1 ({k}_i)$, $\boldsymbol{\Upsilon}_i^2 ({k}_i): \mathbb{U} \mapsto \mathbb{R}^{\mathrm{U}_i \times \mathrm{U}_i}$ are three constant matrices, which are discretized form $\boldsymbol{\Phi}_i({k}_i, {\tau}_i)$, $\boldsymbol{\Phi}^{\prime}_i({k}_i, {\tau}_i)=\partial_{\tau_i} \boldsymbol{\Phi}({k}_i, {\tau}_i)$ and $\boldsymbol{\Phi}^{\prime\prime}_i({k}_i, {\tau}_i)=\partial_{\tau_i}^2 \boldsymbol{\Phi}({k}_i, {\tau}_i)$. After $d$ steps of the $m$-VCF, the NAEs of $\mathbf{R}^d (\mathbf{Q}^d (\boldsymbol{k}_1^d), \boldsymbol{\omega}_1^d) = \mathbf{0} \in \mathbb{R}^{n \mathrm{U}_1^d \times 1}$ with respect to a set of constant coefficients $\mathbf{Q}^d (\boldsymbol{k}_1^d): \mathbb{U}^d \mapsto \mathbb{R}^{n \mathrm{U}_1^d \times 1}$ and base frequencies $\boldsymbol{\omega}_1^d \in \mathbb{R}^{d \times 1}$ are also obtained like Eq. \eqref{EqA2}. 
In this Section, the notation $\square_q^p$ is defined as $\square_q^p = [\square_q; \dots; \square_p]$, and the value $\mathrm{U}_1^{i-1}$ is equal to $\prod_{j=1}^{i-1} \mathrm{U}_j$, with $\mathrm{U}_1^{0}:=1$.

The NAEs of Eq. \eqref{EqA2} of two types of methods both have $n\mathrm{U} = n\mathrm{U}_1^d$ unknown coefficients of $\mathbf{Q}^d$. Both frameworks can be used to construct NAEs based on HB, CO or FD methods. To the best of the authors' knowledge, the VCHB constructed by $m$-VCF can be regarded as the most efficient discretization method for solving quasi-periodic solutions at present.

\section{Construction of $k$-parameters}\label{appB}

In this work, an iterative selection technique is proposed to design the $k$-parameters $\tilde{\boldsymbol{k}}$ in Sub-section \ref{sec3.1}. Assume that the absolutes of $k_j$ collect a vector of $K_j = [0; \boldsymbol{k}_j] \in \mathbb{R}^{(\mathrm{L}_j + 1) \times 1}$, $j = 2, \dots, d$, with $\boldsymbol{k}_j > 0$. Given the parity of the cos and sin functions, a matrix $\tilde{K}^d$ is defined to select $\tilde{\boldsymbol{k}}$ as follows:

a), if $d = 2$, the matrix $\tilde{K}^2 = [\tilde{K}^{2, r}; \tilde{K}^{2, 0}] \in \mathbb{R}^{(\tilde{\mathrm{L}}^2 + 1) \times 1}$ has two parts: namely $\tilde{K}^{2, r} = \boldsymbol{k}_2$ and $\tilde{K}^{2, 0} = 0$, here $\tilde{\mathrm{L}}^2 = \mathrm{L}_2$;

b), if $d > 2$, the matrix $\tilde{K}^d \in \mathbb{R}^{(\tilde{\mathrm{L}}^d + 1) \times (d - 1)}$ is iterated by the matrix $\tilde{K}^{d - 1} \in \mathbb{R}^{(\tilde{\mathrm{L}}^{d - 1} + 1) \times (d - 2)}$, which has two parts:
\begin{equation}
	\tilde{K}^d = \left[ \tilde{K}^{d, g}; \tilde{K}^{d, b} \right],
	\label{EqB1}
\end{equation}
with
\begin{equation}
	\begin{aligned}
		K^{d, g} &= \left[ \boldsymbol{e}_{\mathrm{L}_d} \otimes \tilde{K}^{d - 1, r}, -\boldsymbol{k}_d \otimes \boldsymbol{e}_{\tilde{\mathrm{L}}^{d - 1}} \right]; \\
		\tilde{K}^{d, b} &= \left[ \boldsymbol{e}_{\mathrm{L}_d + 1} \otimes \tilde{K}^{d - 1}, \left[ \boldsymbol{k}_d; 0 \right] \otimes \boldsymbol{e}_{\tilde{\mathrm{L}}^{d - 1} + 1} \right],
	\end{aligned}
	\label{EqB2}
\end{equation}

\begin{figure}[htb] 
	\centering 
	\includegraphics[width=0.95\textwidth]{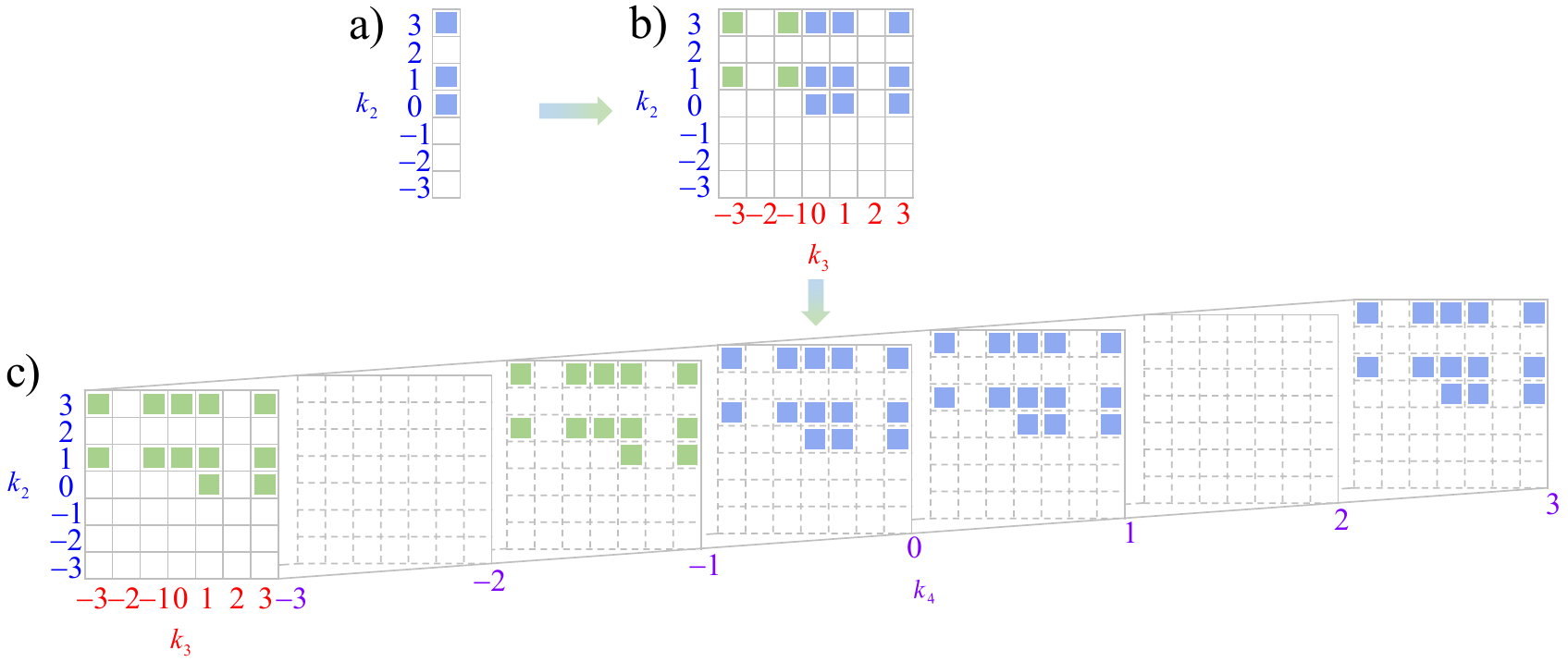} 
	\caption{\small The selected set $\tilde{K}^d$ of harmonic order of $\boldsymbol{\mathcal{H}}(\tilde{\boldsymbol{k}}, \tilde{\boldsymbol{\varphi}})$: a) for $d = 2$; b) for $d = 3$; c) for $d = 4$.} 
	\label{fig: harmonic order}
\end{figure}

where $\tilde{K}^{d-1} = [\tilde{K}^{d-1, r}; \tilde{K}^{d-1, 0}]$ with $\tilde{K}^{d-1, 0}$ and $\tilde{K}^{d-1, r}$ being the remaining harmonic orders. $\boldsymbol{e}_{(\cdot)}$ is all-ones column vector with size of $(\cdot)$. The value $\tilde{\mathrm{L}}^k$ is defined as $[( \prod_{j=2}^k U_k ) - 1]/2$ with $\mathrm{U}_k = 2\mathrm{L}_k + 1$. $\tilde{\boldsymbol{k}}_l$ and $\tilde{\boldsymbol{k}}_0$ are the $l$-th and last row of the matrix $\tilde{K}^d$, respectively. Fig. \ref{fig: harmonic order} provides an example of selection for $\tilde{K}^d$ with $d = 2, 3, 4$, where $K_j$ is defined as $[0, 1, 3]^\mathrm{T}$. The green blocks and blue blocks represent the harmonic order $K^{d, g}$ and $\tilde{K}^{d, b}$ in Eq. \eqref{EqB1}, respectively. Noted that the selection of harmonic order for $d = 2$ and $d = 3$ are same as that of periodic solution in HB [43] and quasi-periodic solution with 2 base frequencies in MHB [30, 39, 66].

\section{One-parameter continuation}\label{appC}

Using $p$ as the continuation parameter, the $(j+1)$-th solution is firstly predicted by the tangent prediction at the known $j$-th solution:
\begin{equation}
	\left[
	\begin{array}{cc}
		\partial_{\boldsymbol{\chi}} \boldsymbol{\mathcal{R}} \big|_{(\boldsymbol{\chi}_j, p_j)} & \partial_p \boldsymbol{\mathcal{R}} \big|_{(\boldsymbol{\chi}_j, p_j)} \\
		\Delta \boldsymbol{\chi}_{j-1}^\top & \Delta p_{j-1}
	\end{array}
	\right]
	\left[
	\begin{array}{c}
		\Delta \boldsymbol{\chi}_j \\
		\Delta p_j
	\end{array}
	\right]
	=
	\left[
	\begin{array}{c}
		0 \\
		1
	\end{array}
	\right]
	\label{EqC1}
\end{equation}
where $\partial_{\boldsymbol{\chi}} \boldsymbol{\mathcal{R}}$, $\partial_p \boldsymbol{\mathcal{R}}$ stands for the derivative of $\boldsymbol{\mathcal{R}}$ with respect to $\boldsymbol{\chi}$, $p$, respectively. With the normalization of $[\Delta \boldsymbol{\chi}_j; \Delta p_j] / |\Delta p_j|$, the prediction of the $(j+1)$-th solution is
\begin{equation}
	\boldsymbol{\chi}_{j+1}^0 = \boldsymbol{\chi}_j + s_j \Delta \boldsymbol{\chi}_j, \quad p_{j+1}^0 = p_j + s_j \Delta p_j
	\label{EqC2}
\end{equation}
where $s_j$ is the step of continuation. Then the solution is corrected by the orthogonal corrections
\begin{equation}
	\left[
	\begin{array}{cc}
		\partial_{\boldsymbol{\chi}} \boldsymbol{\mathcal{R}} \big|_{(\boldsymbol{\chi}_{j+1}^k, p_{j+1}^k)} & \partial_p \boldsymbol{\mathcal{R}} \big|_{(\boldsymbol{\chi}_{j+1}^k, p_{j+1}^k)} \\
		\Delta \boldsymbol{\chi}_j^\top & \Delta p_j
	\end{array}
	\right]
	\left[
	\begin{array}{c}
		\delta \boldsymbol{\chi}_{j+1}^{k+1} \\
		\delta p_{j+1}^{k+1}
	\end{array}
	\right]
	=
	\left[
	\begin{array}{c}
		-\boldsymbol{\mathcal{R}} \big|_{(\boldsymbol{\chi}_{j+1}^k, p_{j+1}^k)} \\
		0
	\end{array}
	\right]
	\label{EqC3}
\end{equation}
\[
\boldsymbol{\chi}_{j+1}^{k+1} = \boldsymbol{\chi}_{j+1}^k + \delta \boldsymbol{\chi}_{j+1}^{k+1}, \quad p_{j+1}^{k+1} = p_{j+1}^k + \delta p_{j+1}^{k+1}
\]
$(\boldsymbol{\chi}_{j+1}^{k+1}, p_{j+1}^{k+1})$ is considered as a solution of Eq. \eqref{EqC3} until $\|\boldsymbol{\mathcal{R}}\|_2 < \varepsilon$ is satisfied, where $\varepsilon$ is a user-defined accuracy.

\section{Alternating Frequency-Time method}\label{appD}

In the subspace $\tilde{\boldsymbol{\varphi}} \in \bar{\mathbb{T}}^{d-1}$, the time samples $\tilde{\boldsymbol{\varphi}}_s$, $s = 1, \ldots, \tilde{\mathrm{S}}$ are defined as the $s$-th column of the matrix $\mathbf{T} = [\mathbf{T}_2; \cdots; \mathbf{T}_d] \in \mathbb{R}^{(d-1) \times \tilde{\mathrm{S}}}$, where
\begin{equation}
	\mathbf{T}_j = \boldsymbol{e}_{\mathrm{S}_{j+1}^d}^{\mathrm{T}} \otimes \bar{\boldsymbol{\varphi}}_j \otimes \boldsymbol{e}_{\mathrm{S}_2^{j-1}}^{\mathrm{T}},
	\label{EqD1}
\end{equation}
Here, $\bar{\boldsymbol{\varphi}}_j \in \mathbb{R}^{1 \times \mathrm{S}_j}$ is the vector of time samples with respect to $\varphi_j$, where $\bar{\boldsymbol{\varphi}}_j = 2\pi \left[ 0, \dots, \mathrm{S}_j - 1 \right]/\mathrm{S}_j,\, j = 2, \dots, d$. $\boldsymbol{e}_{(\cdot)}$ is also an all-ones column vector with size of $(\cdot)$. $(\cdot)^\mathrm{T}$ denotes the transpose operation. $\bar{\mathbf{z}}(\varphi_1, \boldsymbol{\tilde{\varphi}})$ is constructed by $\bar{z}_{i,s} = \bar{z}_i(\varphi_1, \boldsymbol{\tilde{\varphi}}_s) \in \mathbb{R}^{1 \times 1}$ with the rule of $\bar{\mathbf{z}} = \left[ \bar{\mathbf{z}}_1; \dots; \bar{\mathbf{z}}_{2n} \right] \quad \text{and} \quad \bar{\mathbf{z}}_i = [ \bar{z}_{i,1}; \dots; \bar{z}_{i,\tilde{\mathrm{S}}} ],\, i = 1, \dots, 2n$.

The matrices $\boldsymbol{\Gamma} \in \mathbb{R}^{\tilde{\mathrm{S}} \times \tilde{\mathrm{U}}}$ and $\boldsymbol{\Gamma}^{-1} \in \mathbb{R}^{\tilde{\mathrm{U}} \times \tilde{\mathrm{S}}}$ are the constant matrices of the $d$-$1$ dimensional inverse discrete Fourier transforms ($i$-DFT$^{d-1}$) and the $d$-$1$ dimensional discrete Fourier transforms (DFT$^{d-1}$).
\begin{equation}
	\boldsymbol{\Gamma} = \begin{bmatrix}
		1 & \cos(\tilde{\boldsymbol{k}}_1 \tilde{\boldsymbol{\varphi}}_1) & \sin(\tilde{\boldsymbol{k}}_1 \tilde{\boldsymbol{\varphi}}_1) & \dots & \cos(\tilde{\boldsymbol{k}}_{\tilde{\mathrm{L}}^d} \tilde{\boldsymbol{\varphi}}_1) & \sin(\tilde{\boldsymbol{k}}_{\tilde{\mathrm{L}}^d} \tilde{\boldsymbol{\varphi}}_1) \\
		\vdots & \vdots & \vdots & \ddots & \vdots & \vdots \\
		1 & \cos(\tilde{\boldsymbol{k}}_1 \tilde{\boldsymbol{\varphi}}_{\tilde{\mathrm{S}}}) & \sin(\tilde{\boldsymbol{k}}_1 \tilde{\boldsymbol{\varphi}}_{\tilde{\mathrm{S}}}) & \dots & \cos(\tilde{\boldsymbol{k}}_{\tilde{\mathrm{L}}^d} \tilde{\boldsymbol{\varphi}}_{\tilde{\mathrm{S}}}) & \sin(\tilde{\boldsymbol{k}}_{\tilde{\mathrm{L}}^d} \tilde{\boldsymbol{\varphi}}_{\tilde{\mathrm{S}}})
	\end{bmatrix},
	\label{EqD2}
\end{equation}
and
\begin{equation}
	\boldsymbol{\Gamma}^{-1} = \frac{1}{\tilde{\mathrm{S}}} \begin{bmatrix}
		1 & \dots & 1 \\
		2\cos(\tilde{\boldsymbol{k}}_1 \tilde{\boldsymbol{\varphi}}_1) & \dots & 2\cos(\tilde{\boldsymbol{k}}_1 \tilde{\boldsymbol{\varphi}}_{\tilde{\mathrm{S}}}) \\
		2\sin(\tilde{\boldsymbol{k}}_1 \tilde{\boldsymbol{\varphi}}_1) & \dots & 2\sin(\tilde{\boldsymbol{k}}_1 \tilde{\boldsymbol{\varphi}}_{\tilde{\mathrm{S}}}) \\
		\vdots & \ddots & \vdots \\
		2\cos(\tilde{\boldsymbol{k}}_{\tilde{\mathrm{L}}^d} \tilde{\boldsymbol{\varphi}}_1) & \dots & 2\cos(\tilde{\boldsymbol{k}}_{\tilde{\mathrm{L}}^d} \tilde{\boldsymbol{\varphi}}_{\tilde{\mathrm{S}}}) \\
		2\sin(\tilde{\boldsymbol{k}}_{\tilde{\mathrm{L}}^d} \tilde{\boldsymbol{\varphi}}_1) & \dots & 2\sin(\tilde{\boldsymbol{k}}_{\tilde{\mathrm{L}}^d} \tilde{\boldsymbol{\varphi}}_{\tilde{\mathrm{S}}})
	\end{bmatrix}.
	\label{EqD3}
\end{equation}

\section{Newmark integration}\label{appE}

Drop the matrix $\boldsymbol{\Theta}$, Eq. (1) at $\tilde{\boldsymbol{\varphi}} = \tilde{\boldsymbol{\varphi}}_s$, $s = 1, \dots, \tilde{\mathrm{S}}$ can be rewritten as
\begin{equation}
	\begin{aligned}
			&\omega_1^2 \mathbf{M} \mathbf{q}^{\prime\prime}(\varphi_1, \tilde{\boldsymbol{\varphi}}_s) + \omega_1 \mathbf{D} \mathbf{q}^\prime(\varphi_1, \tilde{\boldsymbol{\varphi}}_s) + \mathbf{K} \mathbf{q}(\varphi_1, \tilde{\boldsymbol{\varphi}}_s) \\& + \mathbf{f}_{nl}(\mathbf{q}, \omega_1 \mathbf{q}^\prime) - \mathbf{e}(\varphi_1, \hat{\boldsymbol{\varphi}}_s + \hat{\rho} \varphi_1) = \mathbf{0},
	\end{aligned}
	\label{EqE1}
\end{equation}
where $(\cdot)^\prime$ and $(\cdot)^{\prime\prime}$ denote partial derivatives with respect to hyper-time variable $\varphi_1$. The arguments $\omega_i t$, $i = 1, \dots, e$ in $\mathbf{e}(\boldsymbol{\Omega} \times t)$ are rewritten as $\omega_1 t = \varphi_1$ or $\omega_i t = \varphi_i + \rho_i \varphi_1$, $i = 2, \dots, e$.

For the simplicity of the formulas, $\mathbf{q}(\varphi_1, \tilde{\boldsymbol{\varphi}}_s)$, $\mathbf{q}^\prime(\varphi_1, \tilde{\boldsymbol{\varphi}}_s)$, $\mathbf{q}^{\prime\prime}(\varphi_1, \tilde{\boldsymbol{\varphi}}_s)$ and $\mathbf{e}(\varphi_1, \hat{\boldsymbol{\varphi}}_s + \hat{\rho} \varphi_1)$ are denoted as $\mathbf{q}$, $\mathbf{u}$, $\mathbf{a}$ and $\mathbf{e}$, respectively. And the above system is rewritten as
\begin{equation}
	\omega_1^2 \mathbf{M} \mathbf{a} + \omega_1 \mathbf{D} \mathbf{u} + \mathbf{K} \mathbf{q} + \mathbf{f}_{nl}(\mathbf{q}, \omega_1 \mathbf{u}) - \mathbf{e} = \mathbf{0},
	\label{EqE2}
\end{equation}
Assume there are $\mathrm{S}_1 + 1$ equidistant hyper-time instant $\{ \varphi_{1,k} \}$ in $\varphi_1 = [0, 2\pi]$
\begin{equation}
	\varphi_{1,k} = k \Delta \varphi_1, \quad k = 0, \dots, \mathrm{S}_1,
	\label{EqE3}
\end{equation}
with step of $\Delta \varphi_1 = 2\pi/\mathrm{S}_1$. The velocity $\mathbf{u}_{k+1} = \mathbf{q}'(\varphi_{1,k+1}, \tilde{\boldsymbol{\varphi}}_s)$ and acceleration $\mathbf{a}_{k+1} = \mathbf{q}''(\varphi_{1,k+1}, \tilde{\boldsymbol{\varphi}}_s)$ at $\varphi_{1,k+1}$ are computed by
\begin{equation}
	\begin{aligned}
		\mathbf{u}_{k+1} &= \mathbf{u}_k + \frac{\mathbf{a}_{k+1} + \mathbf{a}_k}{2} \Delta \varphi_1; \\
		\mathbf{q}_{k+1} &= \mathbf{q}_k + \frac{\mathbf{u}_{k+1} + \mathbf{u}_k}{2} \Delta \varphi_1,
	\end{aligned}
	\label{EqE4}
\end{equation}
From this, it follows
\begin{equation}
	\begin{aligned}
		\mathbf{a}_{k+1} &= \frac{4}{\Delta \varphi_1^2} \left( \mathbf{q}_{k+1} - \mathbf{q}_k \right) - \frac{4}{\Delta \varphi_1} \mathbf{u}_k - \mathbf{a}_k; \\
		\mathbf{u}_{k+1} &= \frac{2}{\Delta \varphi_1} \left( \mathbf{q}_{k+1} - \mathbf{q}_k \right) - \mathbf{u}_k,
	\end{aligned}
	\label{EqE5}
\end{equation}
By substituting Eq. \eqref{EqE5} into Eq. \eqref{EqE2}, it is rewritten as
\begin{equation}
	\mathbf{S} \mathbf{q}_{k+1} + \mathbf{f}_{nl} \left( \mathbf{q}_{k+1}, \omega_1 \mathbf{u}_{k+1} (\mathbf{q}_{k+1}) \right) = \mathbf{b},
	\label{EqE6}
\end{equation}
with
\begin{equation}
	\begin{aligned}
		\mathbf{S} &= \omega_1^2 \frac{4}{\Delta \varphi_1^2} \mathbf{M} + \omega_1 \frac{2}{\Delta \varphi_1} \mathbf{D} + \mathbf{K}; \\
		\mathbf{b} &= \mathbf{e}_{k+1} + \omega_1^2 \mathbf{M} \left( \frac{4}{\Delta \varphi_1^2} \mathbf{q}_k + \frac{4}{\Delta \varphi_1} \mathbf{u}_k + \mathbf{a}_k \right) + \omega_1 \mathbf{D} \left( \frac{2}{\Delta \varphi_1} \mathbf{q}_k + \mathbf{u}_k \right),
	\end{aligned}
	\label{EqE7}
\end{equation}
Using Newton iteration for the Eq. \eqref{EqE6} and its adjoint system, $\overline{\mathbf{z}}_{k+1} = \left[ \mathbf{q}_{k+1}, \mathbf{u}_{k+1} \right]$, $\partial_{\overline{\mathbf{z}}_0} \overline{\mathbf{z}}_{k+1}$ and $\partial_{\omega_1} \overline{\mathbf{z}}_{k+1}$ are computed until $k = \mathrm{S}_1$, i.e., $\overline{\mathbf{z}}(2\pi, \tilde{\boldsymbol{\varphi}}_s) \in \mathbb{R}^{2n \times 1}$, $\partial_{\overline{\mathbf{z}}(0, \tilde{\boldsymbol{\varphi}}_s)} \overline{\mathbf{z}}(2\pi, \tilde{\boldsymbol{\varphi}}_s) \in \mathbb{R}^{2n \times 2n}$ and $\partial_{\omega_i} \overline{\mathbf{z}}(2\pi, \tilde{\boldsymbol{\varphi}}_s) \in \mathbb{R}^{2n \times 1},\,i=1,2,\cdots,e$. For example, $\partial_{\overline{\mathbf{z}}(0, \tilde{\boldsymbol{\varphi}}_s)} \overline{\mathbf{z}}(2\pi, \tilde{\boldsymbol{\varphi}}_s)$ is governed by:

\begin{equation}
	\partial_{\overline{\mathbf{z}}(0, \tilde{\boldsymbol{\varphi}}_s)} \overline{\mathbf{z}}(2\pi, \tilde{\boldsymbol{\varphi}}_s) = \begin{bmatrix}
		\frac{\partial \mathbf{q}(2\pi, \tilde{\boldsymbol{\varphi}}_s)}{\partial \mathbf{q}(0, \tilde{\boldsymbol{\varphi}}_s)} & \frac{\partial \mathbf{q}(2\pi, \tilde{\boldsymbol{\varphi}}_s)}{\partial \mathbf{u}(0, \tilde{\boldsymbol{\varphi}}_s)} \\
		\frac{\partial \mathbf{u}(2\pi, \tilde{\boldsymbol{\varphi}}_s)}{\partial \mathbf{q}(0, \tilde{\boldsymbol{\varphi}}_s)} & \frac{\partial \mathbf{u}(2\pi, \tilde{\boldsymbol{\varphi}}_s)}{\partial \mathbf{u}(0, \tilde{\boldsymbol{\varphi}}_s)}
	\end{bmatrix}
	\label{EqE8}
\end{equation}
And, $\overline{\mathbf{z}}(2\pi, \tilde{\boldsymbol{\varphi}}) \in \mathbb{R}^{2n\tilde{\mathrm{S}} \times 1}$, $\partial_{\overline{\mathbf{z}}(0, \tilde{\boldsymbol{\varphi}})} \overline{\mathbf{z}}(2\pi, \tilde{\boldsymbol{\varphi}}) \in \mathbb{R}^{2n\tilde{\mathrm{S}} \times 2n\tilde{\mathrm{S}}}$ and $\partial_{\omega_i} \overline{\mathbf{z}}(2\pi, \tilde{\boldsymbol{\varphi}}) \in \mathbb{R}^{2n\tilde{\mathrm{S}} \times 1}$ are obtained by
\begin{equation}
	\begin{aligned}
		\overline{\mathbf{z}}(2\pi, \tilde{\boldsymbol{\varphi}}) &= \sum_{s=1}^{\tilde{\mathrm{S}}} \overline{\mathbf{z}}(2\pi, \tilde{\boldsymbol{\varphi}}_s) \otimes \mathbf{o}_s; \\
		\partial_{\omega_i} \overline{\mathbf{z}}(2\pi, \tilde{\boldsymbol{\varphi}}) &= \sum_{s=1}^{\tilde{\mathrm{S}}} \partial_{\omega_1} \overline{\mathbf{z}}(2\pi, \tilde{\boldsymbol{\varphi}}_s) \otimes \mathbf{o}_s; i=1,2,\cdots,e\\
		\partial_{\overline{\mathbf{z}}(0, \tilde{\boldsymbol{\varphi}})} \overline{\mathbf{z}}(2\pi, \tilde{\boldsymbol{\varphi}}) &= \sum_{s=1}^{\tilde{\mathrm{S}}} \partial_{\mathbf{z}(0, \tilde{\boldsymbol{\varphi}}_s)} \overline{\mathbf{z}}(2\pi, \tilde{\boldsymbol{\varphi}}_s) \otimes \mathbf{o}_{s,s},
	\end{aligned}
	\label{EqE9}
\end{equation}
Here, $\mathbf{o}_s \in \mathbb{R}^{\tilde{\mathrm{S}}}$ is a column vector with all elements being zero, except that the $m$-th element is equal to $1$. $\mathbf{o}_{s,s}$ is a matrix with all elements being zero, except that the element at $(s,s)$ is equal to $1$. Notedly, the computation of $\overline{\mathbf{z}}(2\pi, \tilde{\boldsymbol{\varphi}}_s)$, $s = 1, \dots, \tilde{\mathrm{S}}$ and its by-production are independent for each $\tilde{\boldsymbol{\varphi}}_s$. So, the procedures of Newmark integration for each $\overline{\mathbf{z}}(2\pi, \tilde{\boldsymbol{\varphi}}_s)$ can be computed by parallelized computation.

\section{Stability analysis based on Lyapunov exponents}\label{appF}

Now, back to the original ODEs in Eq. \eqref{1ODEs}, assuming that $\mathbf{x}_s(t)$ is its stationary solution, the stability of this solution is assessed by the perturbation method. Adding $\Delta \mathbf{x}(0)$ to $\mathbf{x}_s(0)$ as an initial perturbation, i.e., $\mathbf{x}(0) = \mathbf{x}_s(0) + \Delta \mathbf{x}(0)$, and submitting it into Eq. \eqref{1ODEs} yields the perturbation system:
\begin{equation}
	\Delta \dot{\mathbf{x}}(t) \approx \mathbf{J}(t) \Delta \mathbf{x}(t), \quad \mathbf{J}(t) = \left. \frac{\partial \mathbf{f}(\mathbf{x}, t)}{\partial \mathbf{x}(t)} \right|_{\mathbf{x}_s(t)}
	\label{EqF1}
\end{equation}
where $\Delta \mathbf{x}(t)$ represents the evolution of perturbation at time $t$ and equals to $\boldsymbol{\Psi}(t, 0) \Delta \mathbf{x}(0)$ with $\boldsymbol{\Psi}(t, 0) \in \mathbb{R}^{2n \times 2n}$ being the transition matrix from $\Delta \mathbf{x}(0)$ to $\Delta \mathbf{x}(t)$. And $\boldsymbol{\Psi}(t, 0)$ is governed by $\dot{\boldsymbol{\Psi}}(t) \approx \mathbf{J}(t) \boldsymbol{\Psi}(t)$ with initial value $\boldsymbol{\Psi}(0) = \mathbf{I}_{2n}$. The Lyapunov exponent of the asymptotic behavior $t \to +\infty$ of perturbations $\Delta \mathbf{x}(t)$ is usually used to evaluate the stability of quasi-periodic solution \cite{R54,R55}. In this work, the $h$-th order Lyapunov exponents are used for this task based on the discrete Gram-Schmidt orthonormalization:
\begin{equation}
	\begin{aligned}
		\sigma^{(h)}(t_i) &= \lim_{i \to \infty} \frac{1}{t_i} \ln V_h(t_i), \quad h = 1, 2, \dots, 2n; \\
		\sigma_h(t_i) &= \sigma^{(h)}(t_i) - \sigma^{(h-1)}(t_i),
	\end{aligned}
	\label{EqF2}
\end{equation}
where $\sigma^{(h)}(t_i)$ is the $h$-th order Lyapunov exponent, $\sigma_h(t_i)$ is the $h$-th first-order Lyapunov exponent. Here, $t_i$ is equal to $i \Delta t$. If all first-order Lyapunov exponents is not larger than $0$, the quasi-periodic solution is stable; otherwise, it is unstable. Notedly, $V_h(t_i)$ is the volume of $h$-dimensional parallelepiped of the transition matrix $\boldsymbol{\Psi}(t_i, 0)$. Here, $\boldsymbol{\Psi}(t_i, 0)$ is calculated by $\boldsymbol{\Psi}(t_i, t_{i-1}) \hat{\boldsymbol{\Psi}}(t_{i-1}, 0)$, where $\hat{\boldsymbol{\Psi}}(t_{i-1}, 0)$ is re-orthonormalized and re-normalized by the transition matrix $\boldsymbol{\Psi}(t_{i-1}, 0)$ based on the discrete Gram-Schmidt orthonormalization.

\begin{figure}[htb] 
	\centering 
	\includegraphics[width=0.8\textwidth]{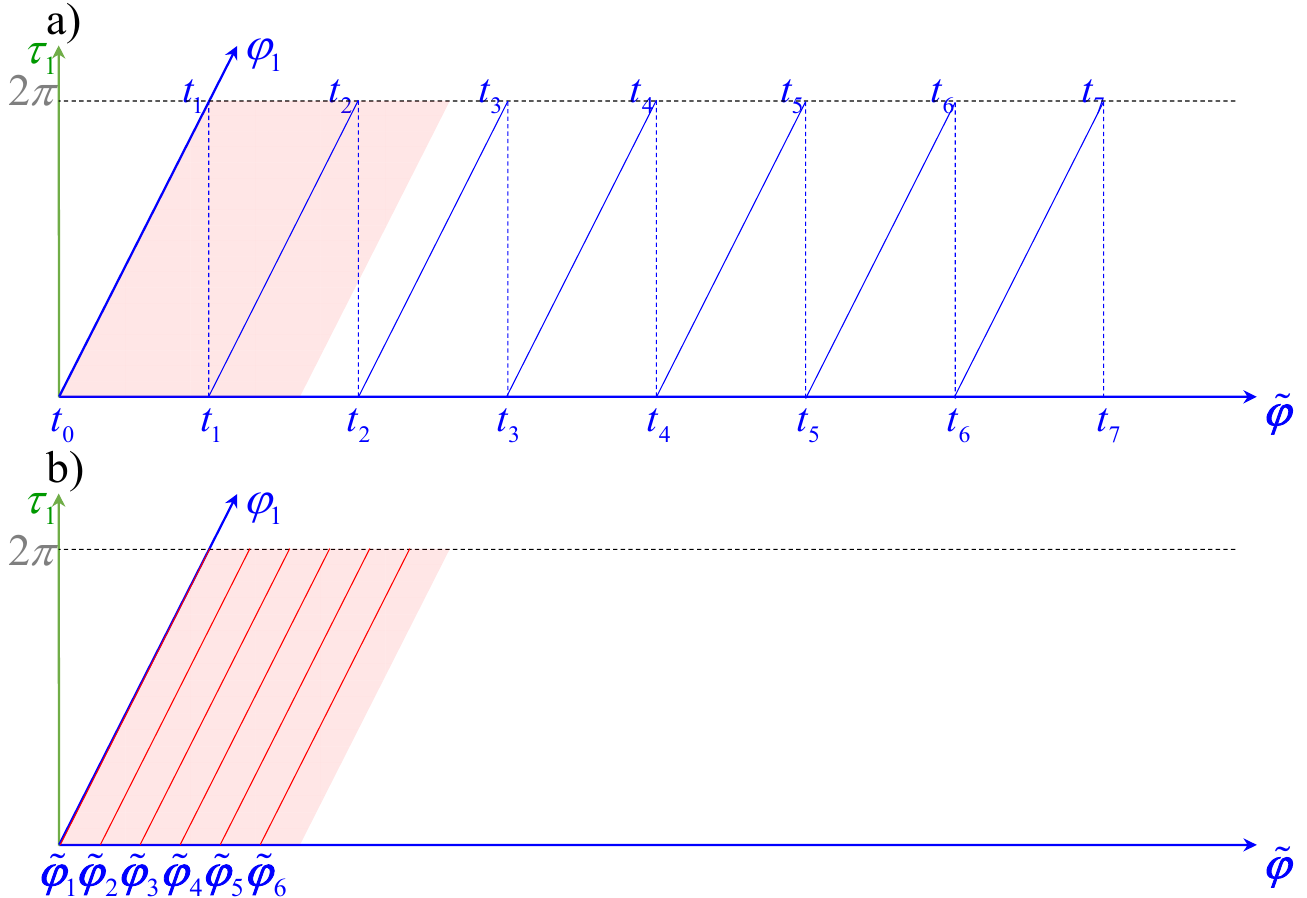} 
	\caption{\small The stability of quasi-periodic solution in term of Lyapunov exponents: (a) The discrete Gram-Schmidt orthonormalization with $\Delta t = T_1$; (b) The by-product of FSE-Shooting for the interpolation of transition matrix $\boldsymbol{\Psi}(t_i, t_{i-1})$.} 
	\label{fig: Lyapunov exponents}
\end{figure}

In this work, the time interval $\Delta t$ is defined as $\Delta t = T_1 = 2\pi/\omega_1$. By introduction the variant of the hyper-time domain of FSE-Shooting, $t_{i-1}$ and $t_i$ can be represented as $\boldsymbol{\varphi}_{i-1} = [2i\pi - 2\pi, 2(i-1)\pi \boldsymbol{\rho}]$ and $\boldsymbol{\varphi}_i = [2i\pi, 2i\pi \boldsymbol{\rho}]$, respectively. Given the periodicity, they can be rewritten as $\boldsymbol{\varphi}_{i-1} = [0, 2(i-1)\pi \boldsymbol{\rho}]$ and $\boldsymbol{\varphi}_i = [2\pi, 2i\pi \boldsymbol{\rho}]$ in Fig. \ref{fig: Lyapunov exponents}a. It can be seen that the transition matrix $\boldsymbol{\Psi}$ is a function of $\boldsymbol{\varphi}$, and moreover, $\tilde{\boldsymbol{\varphi}}$. In the process of Newmark integration, the by-products of $\partial_{\overline{\mathbf{z}}(0, \tilde{\boldsymbol{\varphi}}_s)} \overline{\mathbf{z}}(2\pi, \tilde{\boldsymbol{\varphi}}_s)$, $s = 1, \dots, \tilde{\mathrm{S}}$ in Eq. \eqref{EqChainR} are related with the value of function $\boldsymbol{\Psi}(\tilde{\boldsymbol{\varphi}})$ at $\tilde{\boldsymbol{\varphi}} = \tilde{\boldsymbol{\varphi}}_s$:
\begin{equation}
	\boldsymbol{\Psi}(\tilde{\boldsymbol{\varphi}}_s) = \begin{bmatrix}
		\frac{\partial \mathbf{q}(2\pi, \tilde{\boldsymbol{\varphi}}_s)}{\partial \mathbf{q}(0, \tilde{\boldsymbol{\varphi}}_s)} & \frac{\partial \mathbf{q}(2\pi, \tilde{\boldsymbol{\varphi}}_s)}{\partial \dot{\mathbf{q}}(0, \tilde{\boldsymbol{\varphi}}_s)} \\
		\frac{\partial \dot{\mathbf{q}}(2\pi, \tilde{\boldsymbol{\varphi}}_s)}{\partial \mathbf{q}(0, \tilde{\boldsymbol{\varphi}}_s)} & \frac{\partial \dot{\mathbf{q}}(2\pi, \tilde{\boldsymbol{\varphi}}_s)}{\partial \dot{\mathbf{q}}(0, \tilde{\boldsymbol{\varphi}}_s)}
	\end{bmatrix} = \begin{bmatrix}
		\frac{\partial \mathbf{q}(2\pi, \tilde{\boldsymbol{\varphi}}_s)}{\partial \mathbf{q}(0, \tilde{\boldsymbol{\varphi}}_s)} & \frac{\partial \mathbf{q}(2\pi, \tilde{\boldsymbol{\varphi}}_s)}{\omega_1\partial  \mathbf{u}(0, \tilde{\boldsymbol{\varphi}}_s)} \\
		\frac{\omega_1 \partial \mathbf{u}(2\pi, \tilde{\boldsymbol{\varphi}}_s)}{\partial \mathbf{q}(0, \tilde{\boldsymbol{\varphi}}_s)} & \frac{\partial \mathbf{u}(2\pi, \tilde{\boldsymbol{\varphi}}_s)}{\partial \mathbf{u}(0, \tilde{\boldsymbol{\varphi}}_s)}
	\end{bmatrix}.
	\label{EqF3}
\end{equation}
The transition matrix of $\boldsymbol{\Psi}(\tilde{\boldsymbol{\varphi}})$ is a quasi-periodic function of $\tilde{\boldsymbol{\varphi}}$, as shown in Fig. \ref{fig: Lyapunov exponents}b. So, it can be expressed by means of trigonometric polynomials based on $(d-1)$-dimensional FFT. In this case, the transition matrix $\boldsymbol{\Psi}(t_i, t_{i-1})$ in the discrete Gram-Schmidt orthonormalization is interpolated as $\boldsymbol{\Psi}(\tilde{\boldsymbol{\varphi}}^{i-1})$ with $\tilde{\boldsymbol{\varphi}}^{i-1} = 2(i-1)\pi \boldsymbol{\rho}$. In this work, the maximum of $i$ is denoted as $N_{Ly}$.

\section{Application for self-excited system}\label{appG}

The self-excited system can also exhibits quasi-periodic solutions (e.g., via secondary Hopf bifurcations). A self-excited airfoil system with a store \cite{R33} is used to verify the applicability of proposed phase conditions to self-excited systems, whose governing equation is given as follows:
\begin{equation}
	\mathbf{M}_{\mathbf{q}} \ddot{\mathbf{q}} + \mathbf{C}_{\mathbf{q}} \dot{\mathbf{q}} + \mathbf{K}_{\mathbf{q}} \mathbf{q} + \mathbf{K}_{3}(\mathbf{q}) = \mathbf{0}
	\label{EqG1}
\end{equation}
The degrees of freedom $q_1$, $q_2$, and $q_3$ represent the relative magnitudes of the wing plunge displacement, pitch angle, and store pitch angle, respectively. The related matrices and vectors can be found in \cite{R33}. For this self-excited system, as $Q$ evolves, its stable equilibrium point undergoes a Hopf bifurcation and loses stability at $Q = 4.7440$, inducing a periodic solution with an unknown base frequency; this stable periodic solution then undergoes an NS bifurcation (or secondary Hopf bifurcations) and loses stability at $Q = 7.0769$, inducing a quasi-periodic solution with two unknown base frequencies. 

\begin{figure}[htb] 
	\centering 
	\includegraphics[width=0.85\textwidth]{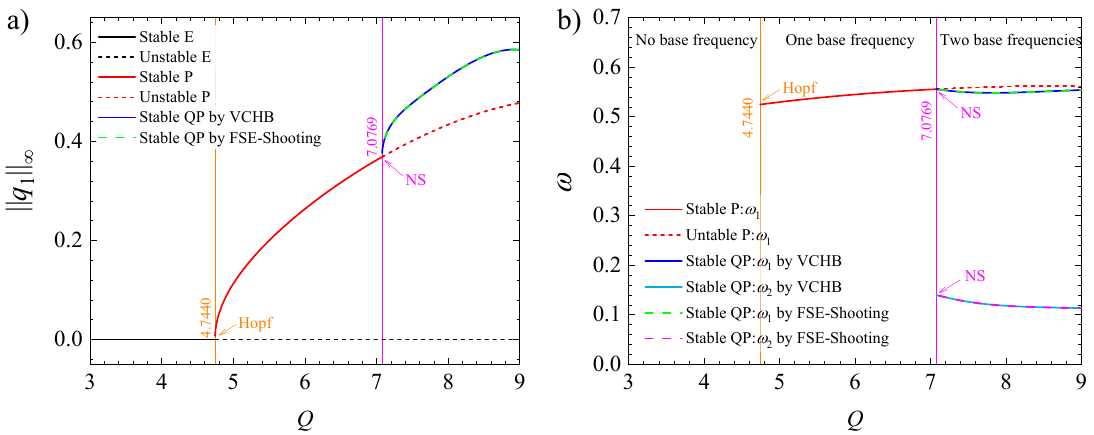} 
	\caption{\small Evolution of  results with respect to the system parameter $Q$ for self-excited airfoil system with a store.: a) the amplitudes of $q_1$; b) base frequencies $\boldsymbol{\omega}$.} 
	\label{fig: self-excited}
\end{figure}

By using Eq. \eqref{EqCont1} in sub-section \ref{sec3.2.1}, Fig. \ref{fig: self-excited} shows the above results computed by VCHB and FSE-Shooting, which is consistent with the descriptions in \cite{R33}. Specifically, Fig. \ref{fig: self-excited}a and b depict the variation of amplitude of $q_1$ and base frequencies of the corresponding solution with the system parameter $Q$, respectively. In detail, when $Q<4.7440$, the system only has stable equilibrium; when $4.7440<Q<7.0769$, the system contains unstable equilibrium and stable periodic solution; subsequently, when $7.0769<Q<9$, the system contains unstable equilibrium point and periodic solution, and stable quasi-periodic solution. Notably, VCHB is faster than FSE-Shooting in this case due to its simpler structure.

\section{Nomenclature}\label{appH}

\begin{table}[ht]
	\resizebox{\linewidth}{!}{
		\centering
		\begin{tabular}{ll}
			\toprule
			Notation & Representation \\
			\midrule
			$n$ & Number of DOFs in $2^\text{nd}$ ODEs \\
			$\mathbf{x}, \mathbf{y}, \mathbf{z}$ & Quasi-periodic state variables in different time domain \\
			$\mathbf{q}$ & Quasi-periodic displacement in $2^\text{nd}$ ODEs \\
			$\mathbf{M}, \mathbf{D}, \mathbf{K}$ & Mass, viscous damping and stiffness matrices \\
			$\boldsymbol{\Theta}$ & Force distribution matrix \\
			$\boldsymbol{\Omega}, \boldsymbol{\omega}$ & Base frequencies of excitation and response \\
			$\boldsymbol{\theta}, \boldsymbol{\tau}, \boldsymbol{\varphi}$ & Variables of $d$-torus in different domain \\
			$\tilde{\boldsymbol{\theta}}, \tilde{\boldsymbol{\tau}}, \tilde{\boldsymbol{\varphi}}$ & Variables of $(d-1)$-torus in different domain \\
			$\boldsymbol{\rho}$ & Ratio of base frequencies \\
			$k_i$ & Discretization parameter associated with the time variable $\tau_i$ or $\varphi_i$ \\
			$\mathrm{U}_i$ & Number of elements of parameter $k_i$ \\
			$\mathrm{S}_i$ & Number of sub-intervals of $\tau_i$ or $\varphi_i$ \\
			$\mathcal{H}$ & Quasi-periodic trigonometric shape functions of $(d-1)$-torus \\
			$\boldsymbol{\Phi}$ & Quasi-periodic shape functions of $d$-torus \\
			$\mathbf{Z}$ & Variable Fourier coefficients of $(d-1)$-torus in domain $\boldsymbol{\varphi}$ \\
			$\hat{\mathbf{Z}}$ & Fourier coefficients of terminal $(d-1)$-torus with phase shift $2\boldsymbol{\rho}\pi$ \\
			$\mathbf{Q}^d$ & Coefficients of $d$-torus in domain $\boldsymbol{\tau}$ \\
			$\mathrm{U}, \tilde{\mathrm{U}}$ & Number of shape functions of $d$-torus or $(d-1)$-torus \\
			$\mathrm{S}, \tilde{\mathrm{S}}$ & Number of time points on $d$-torus or $(d-1)$-torus \\
			$\boldsymbol{\mathcal{R}}$ & Transformed matrix between $\hat{\mathbf{Z}}(2\pi)$ and $\mathbf{Z}(2\pi)$ \\
			$PC_i$ & Phase conditions for base frequency $\omega_i$ \\
			$\nabla_i$ & Derivative matrix in phase conditions \\
			$\boldsymbol{\Gamma}$ & Constant matrix of $i$-DFT$^{d-1}$ \\
			$\boldsymbol{\Gamma}^{-1}$ & Constant matrix of DFT$^{d-1}$ \\
			$\tilde{\Upsilon}^2, \tilde{\Upsilon}^1, \tilde{\Upsilon}^0$ & Discretization matrices in full-discretization methods \\
			$\Upsilon_i^0, \Upsilon_i^1, \Upsilon_i^2$ & Discretization matrices in semi-discretization methods \\
			$K_j$ & Vector collecting the elements of harmonic order of $k_j$ \\
			$\tilde{K}^d$ & Matrix of harmonic order $\tilde{k}$ \\
			$\boldsymbol{\chi}$ & Unknowns in the continuation \\
			$p$ & Continuation parameter \\
			$\mathbf{T}$ & Matrix of sampling points on $(d-1)$-torus \\
			$\mathbf{q}, \mathbf{u}, \mathbf{a}$ & Displacement, velocity and acceleration in Newmark integration \\
			$\Delta \mathbf{x}(t)$ & Perturbation of stationary solution \\
			$\mathbf{J}(t)$ & Jacobian matrix \\
			$\boldsymbol{\Psi}(t, 0)$ & Transition matrix in the perturbation system \\
			$V_h$ & Volume of $h$-dimensional parallelepiped of the transition matrix \\
			$\sigma^{(h)}$ & The $h$-th order Lyapunov exponent \\
			$\sigma_h$ & The $h$-th first-order Lyapunov exponent \\
			$\hat{\boldsymbol{\Psi}}$ & Re-orthonormalized and re-normalized transition matrix \\
			$N_{Ly}$ & Number of periods for computing Lyapunov exponents \\
			\bottomrule
		\end{tabular}
		\label{tab:notation}
	}
\end{table}




\end{document}